# Simultaneous Inference for High-dimensional Linear Models

Xianyang Zhang[*] and Guang Cheng[†]

*Texas A&M University and Purdue University*

February 24, 2016

**Abstract** This paper proposes a bootstrap-assisted procedure to conduct simultaneous inference for high dimensional sparse linear models based on the recent de-sparsifying Lasso estimator (van de Geer et al. 2014). Our procedure allows the dimension of the parameter vector of interest to be exponentially larger than sample size, and it automatically accounts for the dependence within the de-sparsifying Lasso estimator. Moreover, our simultaneous testing method can be naturally coupled with the margin screening (Fan and Lv 2008) to enhance its power in sparse testing with a reduced computational cost, or with the step-down method (Romano and Wolf 2005) to provide a strong control for the family-wise error rate. In theory, we prove that our simultaneous testing procedure asymptotically achieves the pre-specified significance level, and enjoys certain optimality in terms of its power even when the model errors are non-Gaussian. Our general theory is also useful in studying the support recovery problem. To broaden the applicability, we further extend our main results to generalized linear models with convex loss functions. The effectiveness of our methods is demonstrated via simulation studies.

[*]Assistant Professor, Department of Statistics, Texas A&M University, College Station, TX 77843. E-mail: zhangxiany@stat.tamu.edu.
[†]Associate Professor, Department of Statistics, Purdue University, West Lafayette, IN 47906. E-mail: chengg@purdue.edu. Tel: +1 (765) 496-9549. Fax: +1 (765) 494-0558. Research Sponsored by NSF CAREER Award DMS-1151692, DMS-1418042, Simons Fellowship in Mathematics, Office of Naval Research (ONR N00014-15-1-2331) and a grant from Indiana Clinical and Translational Sciences Institute. Guang Cheng was on sabbatical at Princeton while part of this work was carried out; he would like to thank the Princeton ORFE department for its hospitality and support.

# 1  Introduction

High-dimensional statistics has become increasingly popular due to the rapid development of information technologies and their applications in scientific experiments. There is a huge body of work on sparse estimation of high-dimensional models. The statistical properties of these procedures have been extensively studied in the literature; see Bühlmann and van de Geer (2011). The research effort has recently turned to statistical inference such as constructing confidence interval and hypothesis testing for regression coefficients in high-dimensional sparse models. For instance, Wasserman and Roeder (2009), and Meinshausen et al. (2009) proposed significance tests for high-dimensional regression coefficients based on sample splitting. Lockhart et al. (2014) derived a significance test for variables along the Lasso solution path. Lee et al. (2015) proposed an exact post-selection inference procedure based on the idea of polyhedral selection. Another research line for conducting inference is to exploit the idea of low dimensional projection or inverting the Karush-Kuhn-Tucker condition (see e.g. van de Geer et al. 2014, Zhang and Zhang 2014, Javanmard and Montanari 2014, Belloni et al. 2014).

In the high dimensional regime, grouping of variables and exploiting group structure is quite natural (Yuan and Lin 2006, Meier et al. 2008). Leading examples include multifactor analysis-of-variance and additive modeling. From a practical viewpoint, when research interest concerns not only a single variable but rather a group of variables, it seems indispensable to go beyond an approach of inferring individual regression coefficients. The problem of conducting simultaneous inference or inference for groups of variables in high-dimensional models has gained some recent attention. Based on the basis pursuit solution, Meinshausen (2015) proposed an interesting procedure to construct confidence intervals for groups of variables without restrictive assumptions on the design. Mandozzi and Bühlmann (2014) extended the hierarchical testing method in Meinshausen (2008) to the high-dimensional setting, which is able to detect the smallest groups of variables and asymptotically control the family-wise error rate (FWER).

In this paper, we propose a bootstrap-assisted procedure to conduct simultaneous inference in sparse linear models with (possibly) non-Gaussian errors:

$$Y = \mathbf{X}\beta^0 + \epsilon,$$

where $Y$ is a response, $\mathbf{X}$ is a design matrix, and $\beta^0 = (\beta_1^0, \ldots, \beta_p^0)^T$ is a vector of unknown regression coefficients. Specifically, we consider the following simultaneous testing:

$$H_{0,G} : \beta_j^0 = \widetilde{\beta}_j \quad \text{for all } j \in G \subseteq \{1, 2, \ldots, p\}$$

versus the alternative $H_{a,G} : \beta_j^0 \neq \widetilde{\beta}_j$ for some $j \in G$, where $\widetilde{\beta}_j$ with $j \in G$ is a vector of pre-specified values (e.g., by domain experts). We point out that two extreme cases, where $G = \{1, 2, \ldots, p\}$ or $|G|$ is small, have been considered in the recent literature, see e.g., Arias-Castro et al. (2011), Zhong and Chen (2011), Feng et al. (2013), Zhang and Zhang (2014), and van de Geer et al. (2014). However, relatively few attention has been paid to the case where $G$ lies between these two extremes. Our method allows $G$ to be an arbitrary subset of $\{1, 2, \ldots, p\}$, and thus can be applied to a wider range of real problems such as testing the significance of a growing set of genes (associated with certain clinical outcome) conditional on another set of genes, whose size is allowed to grow as well. In comparison with Meinshausen (2015) and Mandozzi and Bühlmann (2014), our method is asymptotically exact and relies on suitable assumptions (Assumptions 2.1 & 2.2) on the design which are different from the highly correlated design considered in Mandozzi and Bühlmann (2014).

Our general framework is built upon the de-sparsifying Lasso estimator, denoted as $\breve{\beta} = (\breve{\beta}_1, \ldots, \breve{\beta}_p)^T$, recently developed in van de Geer et al. (2014) and Zhang and Zhang (2014), whose review is given in Section 2.1. In Section 2.2, a test statistic is proposed as $T_{n,G} := \max_{j \in G} \sqrt{n} |\breve{\beta}_j - \widetilde{\beta}_j|$, whose critical values are obtained via a simple multiplier bootstrap method. Based on the asymptotic linear expansion of $\breve{\beta}$, we show that the proposed multiplier bootstrap method consistently approximates the null limiting distribution of $T_{n,G}$, and thus the testing procedure achieves the pre-specified significance level asymptotically. It is worth mentioning that the proposed bootstrap-assisted procedure is adaptive to the dimension of the component of interest, and it automatically accounts for the dependence within the de-sparsifying Lasso estimators. In theory, we also prove that our testing procedure enjoys certain minimax optimality (Verzelen 2012) in terms of its power even when the model errors are non-Gaussian and the cardinality of $G$ is exponentially larger than sample size.

Moreover, our new methodology is readily applicable to some other important statistical problems in the high-dimensional setting, such as support recovery, testing for sparse signals, and multiple testing. The support recovery procedure is proposed in Section 3.1 as an important by-product of our general theory. It has been shown through simulations that the proposed procedure can be more accurate in recovering signals than Lasso, the stability selection (Meinshausen and Bülmann 2010), and the screen and clean procedure (Wasserman and Roeder 2009). The above bootstrap-assisted test method can also be coupled with the margin screening in Fan and Lv (2008) to enhance the power performance in sparse testing with a reduced computational cost, which is very attractive in the ultra-high dimensional setting. Hence, in Section 3.2 we propose a three-step procedure that first randomly splits the sample into two subsamples, screens out the irrelevant variables based on the first subsample, and finally performs the above maximum-type testing on the reduced model based on the second subsample. Another application is a multiple testing problem: for each $j \in G$

$$H_{0,j} : \beta_j^0 \leq \widetilde{\beta}_j \text{ versus } H_{a,j} : \beta_j^0 > \widetilde{\beta}_j.$$

To obtain a strong control of the FWER, we incorporate the above bootstrap idea into the step-down method (Romano and Wolf 2005) in Section 3.3. As noted in Chernozhukov et al. (2013), this hybrid method is asymptotically non-conservative as compared to the Bonferroni-Holm procedure since the correlation amongst the test statistics has been taken into account.

To broaden the applicability, we further extend our main results to generalized linear models with convex loss functions in Section 4. The usefulness of the above simultaneous inference methods is illustrated via simulation studies in Section 5. The technical details and additional numerical results are included in a supplement file. Some theoretical derivations in this paper rely on an impressive Gaussian approximation (GAR) theory recently developed in Chernozhukov et al. (2013). The application of GAR theory is nontrivial as one needs to verify suitable moment conditions on the leading term of the de-sparsifying Lasso estimator, and quantify the estimation effect as well as the impact of the remainder term in (7) below. We also want to point out that the GAR theory is applied without conditioning on the design as random design is considered throughout the paper. Our results complement Belloni et al. (2014) who establish the validity of uniform confidence band in the high-dimensional least absolute deviation regression.

Finally, we introduce some notation. For a $p \times p$ matrix $B = (b_{ij})_{i,j=1}^{p}$, let $||B||_{\infty} = \max_{1 \leq i,j \leq p} |b_{ij}|$ and $||B||_1 = \max_{1 \leq j \leq p} \sum_{i=1}^{p} |b_{ij}|$. Denote by $||a||_q = (\sum_{i=1}^{p} |a_i|^q)^{1/q}$ and $||a||_{\infty} = \max_{1 \leq j \leq p} |a_j|$ for $a = (a_1, \ldots, a_p)^T \in \mathbb{R}^p$ and $q > 0$. For a set $\mathcal{A}$, denote its cardinality by $|\mathcal{A}|$. Denote by $\lfloor a \rfloor$ the integer part of a positive real number $a$. For two sequences $\{a_n\}$ and $\{b_n\}$, write $a_n \asymp b_n$ if there exist positive constants $c$ and $C$ such that $c \leq \liminf_n (a_n/b_n) \leq \limsup_n (a_n/b_n) \leq C$. Also write $a_n \lesssim b_n$ if $a_n \leq C' b_n$ for some constant $C' > 0$ independent of $n$ (and $p$). The symbol $N_p(\mu, \Sigma)$ is reserved for a $p$-dimensional multivariate normal distribution with mean $\mu$ and covariance matrix $\Sigma$.

## 2 Main Theory

### 2.1 De-sparsifying Lasso estimator

In this section, we review the de-sparsifying (de-biased) estimator proposed in van de Geer et al. (2014), which is essentially the same as the estimator proposed in Zhang and Zhang (2014) but motivated from a different viewpoint, i.e., inverting the Karush-Kuhn-

Tucker (KKT) condition of Lasso. Consider a high-dimensional sparse linear model:

$$Y = \mathbf{X}\beta^0 + \epsilon, \tag{1}$$

with a response $Y = (Y_1, \ldots, Y_n)^T$, an $n \times p$ design matrix $\mathbf{X} := [X_1, \ldots, X_p]$, an error $\epsilon = (\epsilon_1, \ldots, \epsilon_n)^T$ independent of $\mathbf{X}$, and an unknown $p \times 1$ regression vector $\beta^0 = (\beta_1^0, \ldots, \beta_p^0)^T$. The parameter dimension $p$ can be much larger than sample size $n$. Suppose $\mathbf{X}$ has i.i.d rows having mean zero and covariance matrix $\Sigma = (\sigma_{ij})_{i,j=1}^p$ with $\Sigma^{-1} := \Theta = (\theta_{ij})_{i,j=1}^p$. We denote the active set of variables by $\mathcal{S}_0 = \{1 \leq j \leq p : \beta_j^0 \neq 0\}$ and its cardinality by $s_0 = |\mathcal{S}_0|$. The Lasso estimator (Tibshirani, 1996) is defined as

$$\widehat{\beta} = \text{argmin}_{\beta \in \mathbb{R}^p}(||Y - \mathbf{X}\beta||_2^2/n + 2\lambda ||\beta||_1), \tag{2}$$

for some tuning parameter $\lambda > 0$.

The de-sparsifying estimator is obtained by inverting the KKT condition,

$$\breve{\beta} = \widehat{\beta} + \widehat{\Theta}\mathbf{X}^T(Y - \mathbf{X}\widehat{\beta})/n, \tag{3}$$

where $\widehat{\Theta}$ is a suitable approximation for the inverse of the Gram matrix $\widehat{\Sigma} := \mathbf{X}^T\mathbf{X}/n$. In what follows, we consider the approximate inverse $\widehat{\Theta}$ given by Lasso for the nodewise regression on the design matrix $\mathbf{X}$; see Meinshausen and Bühlmann (2006). Let $\mathbf{X}_{-j}$ be the design matrix without the $j$th column. For $j = 1, 2, \ldots, p$, consider

$$\widehat{\gamma}_j := \text{argmin}_{\gamma \in \mathbb{R}^{p-1}}(||X_j - \mathbf{X}_{-j}\gamma||_2^2/n + 2\lambda_j ||\gamma||_1) \tag{4}$$

with $\lambda_j > 0$, where we denote $\widehat{\gamma}_j = \{\widehat{\gamma}_{j,k} : 1 \leq k \leq p, \ k \neq j\}$. Let $\widehat{C} = (\widehat{c}_{i,j})_{i,j=1}^p$ be a $p \times p$ matrix with $\widehat{c}_{i,i} = 1$ and $\widehat{c}_{i,j} = -\widehat{\gamma}_{i,j}$ for $i \neq j$. Let $\widehat{\tau}_j^2 = ||X_j - \mathbf{X}_{-j}\widehat{\gamma}_j||_2^2/n + \lambda_j ||\widehat{\gamma}_j||_1$ and write $\widehat{T}^2 = \text{diag}(\widehat{\tau}_1^2, \ldots, \widehat{\tau}_p^2)$ as a diagonal matrix. Finally, the nodewise Lasso estimator for

$\Theta$ is constructed as $\widehat{\Theta} = \widehat{T}^{-2}\widehat{C}$.

Denote by $\gamma_j = \text{argmin}_{\gamma \in \mathbb{R}^{p-1}} \mathbb{E}||X_j - \mathbf{X}_{-j}\gamma||_2^2$, and define $\eta_j = X_j - \mathbf{X}_{-j}\gamma_j = (\eta_{1,j}, \ldots, \eta_{n,j})^T$. Define $\tau_j^2 = \mathbb{E}||\eta_j||_2^2/n = 1/\theta_{j,j}$ for $j = 1, 2, \ldots, p$. Let $s_j = |\{1 \leq k \leq p : k \neq j, \theta_{jk} \neq 0\}|$. Denote the $j$th row of $\mathbf{X}$ and $\widehat{\Theta}$ by $\widetilde{X}_j = (X_{j1}, \ldots, X_{jp})^T$ and $\widehat{\Theta}_j$, respectively.

ASSUMPTION 2.1. *The design matrix $\mathbf{X}$ has either i.i.d sub-Gaussian rows (i.e., $\sup_{||a||_2 \leq 1} \mathbb{E} \exp\{|\sum_{j=1}^p a_j X_{ij}|^2/C\} \leq 1$ for some large enough positive constant $C$) or i.i.d rows satisfying for some $K_n \geq 1$, $\max_{1 \leq i \leq n, 1 \leq j \leq p} |X_{ij}| \leq K_n$ (strongly bounded case), where $K_n$ is allowed to grow with $n$. In the strongly bounded case, we assume in addition that $\max_j ||\mathbf{X}_{-j}\gamma_j||_\infty \leq K_n$ and $\max_j \mathbb{E}\eta_{1,j}^4 \leq K_n^4$.*

ASSUMPTION 2.2. *The smallest eigenvalue $\Lambda_{\min}^2$ of $\Sigma$ satisfies that $c < \Lambda_{\min}^2$, and $\max_j \Sigma_{j,j} \leq C$, where $c, C$ are some positive constants.*

THEOREM 2.1 (Theorem 2.4 of van de Geer et al. 2014). *Suppose Assumptions 2.1-2.2 hold. Assume in the sub-Gaussian case, it holds that $\max_{1 \leq j \leq p} \sqrt{s_j \log(p)/n} = o(1)$ and in the strongly bounded case, $\max_{1 \leq j \leq p} K_n^2 s_j \sqrt{\log(p)/n} = o(1)$. Then with suitably chosen $\lambda_j \asymp K_0\sqrt{\log(p)/n}$ uniformly for $j = 1, 2, \ldots, p$, where $K_0 = 1$ in the sub-Gaussian case and $K_0 = K_n$ in the strongly bounded case, we have*

$$||\widehat{\Theta}_j - \Theta_j||_1 = O_P(K_0 s_j \sqrt{\log(p)/n}), \quad ||\widehat{\Theta}_j - \Theta_j||_2 = O_P(K_0\sqrt{s_j \log(p)/n}), \quad (5)$$

$$|\widehat{\tau}_j^2 - \tau_j^2| = O_P(K_0\sqrt{s_j \log(p)/n}), \quad (6)$$

*uniformly for $1 \leq j \leq p$. Furthermore, suppose $\{\epsilon_i\}$ are i.i.d with $c' < \sigma_\epsilon^2 < c$ and in the sub-Gaussian case for $\mathbf{X}$, $\mathbb{E}\exp(|\epsilon_i|/C) \leq 1$ for some positive constants $c, c', C > 0$. Assume that $\lambda$ is suitably chosen such that $\lambda \asymp K_0\sqrt{\log(p)/n}$, and $K_0 s_0 \log(p)/\sqrt{n} = o(1)$ and $\max_{1 \leq j \leq p} K_0 s_j \sqrt{\log(p)/n} = o(1)$. Then*

$$\sqrt{n}(\breve{\beta} - \beta^0) = \widehat{\Theta}\mathbf{X}^T\epsilon/\sqrt{n} + \Delta, \quad ||\Delta||_\infty = o_P(1), \quad (7)$$

*where* $\Delta = (\Delta_1, \ldots, \Delta_p)^T = -\sqrt{n}(\widehat{\Theta}\widehat{\Sigma} - I)(\widehat{\beta} - \beta^0)$.

Theorem 2.1 provides an explicit expansion for the de-sparsifying Lasso estimator and states that the remainder term $\Delta$ can be well controlled in the sense that $||\Delta||_\infty = o_P(1)$, which is very useful in the subsequent derivation.

## 2.2 Simultaneous inference procedures

In the high dimensional regime, it is natural to test the hypothesis

$$H_{0,G} : \beta_j^0 = \widetilde{\beta}_j \quad \text{for all } j \in G \subseteq \{1, 2, \ldots, p\}$$

versus the alternative $H_{a,G} : \beta_j^0 \neq \widetilde{\beta}_j$ for some $j \in G$. For example, we consider the sparse testing, i.e., $\widetilde{\beta}_j = 0$, in Section 3.2. Throughout the paper, we allow $|G|$ to grow as fast as $p$, which can be of an exponential order w.r.t. $n$. Hence, our results go beyond the existing ones in van de Geer et al. (2014), Zhang and Zhang (2014), and Javanmard and Montanari (2014), where $|G|$ is fixed. Simultaneous inference has been considered earlier via Bonferroni adjustment, which leads to very conservative procedures. In contrast, our method is asymptotically nonconservative.

In this section, we propose the test statistic

$$\max_{j \in G} \sqrt{n}|\breve{\beta}_j - \widetilde{\beta}_j|$$

with $\breve{\beta}_j$ being the de-sparsifying estimator. As will be seen in Sections 3.2 and 3.3, this test statistic can be naturally coupled with the margin screening (Fan and Lv 2008) to enhance its power in sparse testing and also reduce the computational cost in nodewise Lasso, or with the step-down method (Romano and Wolf 2005) to provide a strong control for the FWER. We next describe a simple multiplier bootstrap method to obtain an accurate critical value.

The asymptotic linear expansion in (7) can be re-written as

$$(\widehat{\Theta}\mathbf{X}^T\epsilon/\sqrt{n})_j = \sum_{i=1}^n \widehat{\Theta}_j^T \widetilde{X}_i \epsilon_i/\sqrt{n} = \sum_{i=1}^n \widehat{\xi}_{ij}/\sqrt{n}, \quad \widehat{\xi}_{ij} = \widehat{\Theta}_j^T \widetilde{X}_i \epsilon_i,$$

where $(a)_j = a_j$ for $a = (a_1, \ldots, a_p)^T$. Generate a sequence of random variables $\{e_i\}_{i=1}^n \overset{i.i.d.}{\sim} N(0,1)$ and define the multiplier bootstrap statistic,

$$W_G = \max_{j \in G} \sum_{i=1}^n \widehat{\Theta}_j^T \widetilde{X}_i \widehat{\sigma}_\epsilon e_i/\sqrt{n},$$

where $\widehat{\sigma}_\epsilon^2$ is a consistent estimator of the error variance $\sigma_\epsilon^2$, e.g., the variance estimator from the scaled Lasso (Sun and Zhang 2012). The bootstrap critical value is given by $c_G(\alpha) = \inf\{t \in \mathbb{R} : P(W_G \leq t | (Y, \mathbf{X})) \geq 1 - \alpha\}$.

The validity of the above bootstrap method requires the following two assumptions.

ASSUMPTION 2.3. *(i) If $\mathbf{X}$ has i.i.d sub-Gaussian rows, assume that $(\log(pn))^7/n \leq C_1 n^{-c_1}$ for some constants $c_1, C_1 > 0$. In this case, suppose $\{\epsilon_i\}$ are i.i.d sub-Gaussian with $c' < \sigma_\epsilon^2 < c$ for $c, c' > 0$. (ii) If $\max_{1 \leq i \leq n, 1 \leq j \leq p} |X_{ij}| \leq K_n$, assume that $\max_{1 \leq j \leq p} s_j K_n^2 (\log(pn))^7/n \leq C_2 n^{-c_2}$ for some constants $c_2, C_2 > 0$. In this case, suppose $\{\epsilon_i\}$ are i.i.d sub-exponential, i.e., $\mathbb{E}\exp(|\epsilon_i|/C') \leq 1$ and $c' < \sigma_\epsilon^2 < c$ for some constants $c, c', C' > 0$.*

ASSUMPTION 2.4. *There exists a sequence of positive numbers $\alpha_n \to +\infty$ such that $\alpha_n/p = o(1)$, $\alpha_n (\log p)^2 \max_j \lambda_j \sqrt{s_j} = o(1)$ and $P(\alpha_n (\log p)^2 |\widehat{\sigma}_\epsilon^2 - \sigma_\epsilon^2| > 1) \to 0$.*

Assumption 2.3 requires that (i) the regressors and errors are both sub-Gaussian or (ii) the regressors are strongly bounded while the errors are sub-exponential. It also imposes suitable restrictions on the growth rate of $p$. We point out that the existence of a factor $(\log(p))^2$ in Assumption 2.4 is due to an application of Theorem 2.1 in Chernozhukov et al. (2014) regarding the comparison for the maxima of two Gaussian random vectors. Assumption 2.4 is a very mild technical condition. For example if $|\widehat{\sigma}_\epsilon^2 - \sigma_\epsilon^2| = O_P(1/\sqrt{n})$ and $\lambda_j \asymp K_0 \sqrt{\log(p)/n}$ uniformly for all $j$, then Assumption 2.4 holds provided that $K_0 (\log p)^{5/2} \max_j \sqrt{s_j/n} = o(1)$.

It is worth noting that the $\sqrt{n}$ convergence rate for $\widehat{\sigma}_\epsilon^2$ can be achieved in the high dimensional setting (e.g. by the scaled Lasso in Sun and Zhang 2011). Recall that $K_0 = 1$ in the sub-Gaussian case and $K_0 = K_n$ in the strongly bounded case. Theorem 2.2 below establishes the validity of the bootstrap procedure for one sided test.

THEOREM 2.2. *Suppose that Assumptions 2.1-2.4 hold. Assume that*

$$\max_{1\leq j\leq p} K_0^2 s_j^2 (\log(pn))^3 (\log(n))^2 / n = o(1), \tag{8}$$

$$K_0^4 s_0^2 (\log(p))^3 / n = o(1), \tag{9}$$

*and $\lambda$ and $\lambda_j$ are suitably chosen such that*

$$\lambda \asymp K_0 \sqrt{\log(p)/n}, \quad \lambda_j \asymp K_0 \sqrt{\log(p)/n} \text{ uniformly for } j. \tag{10}$$

*Then we have for any $G \subseteq \{1, 2, \ldots, p\}$,*

$$\sup_{\alpha \in (0,1)} \left| P\left( \max_{j \in G} \sqrt{n}(\breve{\beta}_j - \beta_j^0) > c_G(\alpha) \right) - \alpha \right| = o(1). \tag{11}$$

The scaling condition (8) is imposed to control the estimation effect caused by replacing $\Theta$ with its nodewise Lasso estimator $\widehat{\Theta}$, while condition (9) is needed to bound the remainder term $\Delta$. One crucial feature of our bootstrap-assisted testing procedure is that it explicitly accounts for the effect of $|G|$ in the sense that the bootstrap critical value $c_G(\alpha)$ depends on $G$. This is in sharp contrast with the extreme value approach in Cai et al. (2014); see Remark 2.4. Hence, our approach is more robust to the change in $|G|$. Since $\max_{j\in G} \sqrt{n}|\breve{\beta}_j - \widetilde{\beta}_j| = \sqrt{n} \max_{j \in G} \max\{\breve{\beta}_j - \widetilde{\beta}_j, \widetilde{\beta}_j - \breve{\beta}_j\}$, similar arguments imply that under $H_{0,G}$,

$$\sup_{\alpha \in (0,1)} \left| P\left( \max_{j \in G} \sqrt{n}|\breve{\beta}_j - \widetilde{\beta}_j| > c_G^*(\alpha) \right) - \alpha \right| = o(1),$$

where $c_G^*(\alpha) = \inf\{t \in \mathbb{R} : P(W_G^* \leq t|(Y,\mathbf{X})) \geq 1 - \alpha\}$ with $W_G^* =$

$\max_{j\in G}|\sum_{i=1}^{n}\widehat{\Theta}_{j}^{T}\widetilde{X}_{i}\widehat{\sigma}_{\epsilon}e_{i}/\sqrt{n}|$. This result is readily applicable to construct simultaneous confidence intervals for $\beta_{j}^{0}$ with $j \in G$.

REMARK 2.1. Alternatively, we can employ Efron's empirical bootstrap to obtain the critical value. For simplicity, we take $G = \{1, 2, \ldots, p\}$. Let $\widehat{h}_{i} = (\widehat{h}_{i1}, \ldots, \widehat{h}_{ip})^{T}$ with $\widehat{h}_{ij} = \widehat{\Theta}_{j}^{T}\widetilde{X}_{i}\widehat{\sigma}_{\epsilon}$. Let $\widehat{h}_{1}^{*}, \ldots, \widehat{h}_{n}^{*}$ be a sample from the empirical distribution based on $\{\widehat{h}_{i}\}_{i=1}^{n}$. Define the empirical bootstrap statistic as $W_{EB}^{*} = \max_{1 \leq j \leq p} |\sum_{i=1}^{n}(\widehat{h}_{ij}^{*} - \sum_{i=1}^{n}\widehat{h}_{ij}/n)/\sqrt{n}|$. The empirical bootstrap critical value is then given by $c_{EB}^{*}(\alpha) = \inf\{t \in \mathbb{R} : P(W_{EB}^{*} \leq t|(Y, \mathbf{X})) \geq 1 - \alpha\}$. Following the arguments in Appendix K of Chernozhukov et al. (2013) and the proof of Theorem 2.2, we can establish the asymptotic equivalence between the empirical bootstrap and the multiplier bootstrap, which theoretically justifies the use of the former. We omit the technical details here to conserve space.

Let $\widehat{\Xi} = (\widehat{\omega}_{ij})_{i,j=1}^{p} = \widehat{\Theta}\widehat{\Sigma}\widehat{\Theta}^{T}\widehat{\sigma}_{\epsilon}^{2}$. We next consider the studentized statistic $\max_{j\in G}\sqrt{n}(\breve{\beta}_{j} - \widetilde{\beta}_{j})/\sqrt{\widehat{\omega}_{jj}}$ for one sided test. In this case, the bootstrap critical value can be obtained via $\bar{c}_{G}(\alpha) = \inf\{t \in \mathbb{R} : P(\bar{W}_{G} \leq t|(Y, \mathbf{X})) \geq 1 - \alpha\}$, where $\bar{W}_{G} = \max_{j\in G}\sum_{i=1}^{n}\widehat{\Theta}_{j}^{T}\widetilde{X}_{i}\widehat{\sigma}_{\epsilon}e_{i}/\sqrt{n\widehat{\omega}_{jj}}$. The following theorem justifies the validity of the bootstrap procedure for the studentized statistic.

THEOREM 2.3. *Under the assumptions in Theorem 2.2, we have for any $G \subseteq \{1, 2, \ldots, p\}$,*

$$\sup_{\alpha \in (0,1)} \left| P\left(\max_{j\in G}\sqrt{n}(\breve{\beta}_{j} - \beta_{j}^{0})/\sqrt{\widehat{\omega}_{jj}} > \bar{c}_{G}(\alpha)\right) - \alpha \right| = o(1). \tag{12}$$

Following the arguments in the proof of Theorem 2.3, it is straightforward to show that under $H_{0,G}$,

$$\sup_{\alpha \in (0,1)} \left| P\left(\max_{j\in G}\sqrt{n}|\breve{\beta}_{j} - \widetilde{\beta}_{j}|/\sqrt{\widehat{\omega}_{jj}} > \bar{c}_{G}^{*}(\alpha)\right) - \alpha \right| = o(1),$$

where $\bar{c}_{G}^{*}(\alpha) = \inf\{t \in \mathbb{R} : P(\bar{W}_{G}^{*} \leq t|(Y, \mathbf{X})) \geq 1 - \alpha\}$ with $\bar{W}_{G}^{*} = \max_{j\in G}|\sum_{i=1}^{n}\widehat{\Theta}_{j}^{T}\widetilde{X}_{i}\widehat{\sigma}_{\epsilon}e_{i}/\sqrt{n\widehat{\omega}_{jj}}|$. Compared to Theorem 2 in Zhang and Zhang (2014), our two-sided testing procedure above is straightforward to implement and asymptotically exact.

In particular, the unknown quantile parameter in (29) of Zhang and Zhang (2014) seems not directly obtainable. Our procedure is also asymptotically nonconservative as compared to the method in Meinshausen (2015).

REMARK 2.2. *Our inferential procedures allow the precision matrix to be sparse. In fact, they work for any estimation method for precision matrix as long as the estimation effect such as $||\widehat{\Theta}^T - \Theta^T||_1$ and $||\widehat{\Theta}\widehat{\Sigma} - I||_\infty$ can be well controlled under suitable assumptions (see e.g. Cai et al. 2011, Javanmard and Montanari 2014).*

REMARK 2.3. *An alternative way to conduct simultaneous inference is based on the sum of squares type statistics, e.g., Chen and Qin (2010), which is expected to have good power against non-sparse alternatives. However, such a type of test statistics may not work well in the current setting due to the accumulation of estimation errors especially when $|G|$ is large, i.e., the error term $\sum_{j \in G} |\Delta_j|$ might be out of control.*

We next turn to the (asymptotic) power analysis of the above procedure. Note that when $|G|$ is fixed, our procedure is known to be $\sqrt{n}$-consistent (implicitly implied by Theorem 2.1). In fact, even when $|G| \to \infty$, our test still enjoys certain optimality in the sense that the separation rate $(\sqrt{2}+\varepsilon_0)\sqrt{\log(|G|)/n}$ for any $\varepsilon_0 > 0$ derived in Theorem 2.4 below is minimax optimal according to Section 3.2 of Verzelen (2012) under suitable assumptions.

Below we focus on the case where $|G| \to \infty$ as $n \to \infty$. Define the separation set

$$\mathcal{U}_G(c_0) = \{\beta = (\beta_1, \ldots, \beta_p)^T : \max_{j \in G} |\beta_j - \widetilde{\beta}_j|/\sqrt{\omega_{jj}} > c_0\sqrt{\log(|G|)/n}\}, \tag{13}$$

where $\omega_{jj} = \sigma_\epsilon^2 \theta_{jj}$. Recall that $\Sigma^{-1} = \Theta = (\theta_{ij})_{i,j=1}^p$. Let $\widetilde{\Theta} = (\widetilde{\theta}_{ij})$ with $\widetilde{\theta}_{ij} = \theta_{ij}/\sqrt{\theta_{ii}\theta_{jj}}$.

ASSUMPTION 2.5. *Assume that $\max_{1 \leq i \neq j \leq p} |\widetilde{\theta}_{ij}| \leq c < 1$ for some constant $c$.*

THEOREM 2.4. *Under the assumptions in Theorem 2.3 and Assumption 2.5, we have for any $\varepsilon_0 > 0$,*

$$\inf_{\beta^0 \in \mathcal{U}_G(\sqrt{2}+\varepsilon_0)} P\left(\max_{j \in G} \sqrt{n}|\breve{\beta}_j - \widetilde{\beta}_j|/\sqrt{\widehat{\omega}_{jj}} > \breve{c}_G^*(\alpha)\right) \to 1. \tag{14}$$

Theorem 2.4 says that the correct rejection of our bootstrap-assisted test can still be triggered even when there exists only one entry of $\beta^0 - \widetilde{\beta}$ with a magnitude being larger than $(\sqrt{2} + \varepsilon_0)\sqrt{\log(|G|)/n}$ in $G$. Hence, our procedure is very sensitive in detecting sparse alternatives. As pointed out by one referee, when $\Sigma$ is an identity matrix, the constant $\sqrt{2}$ turns out to be asymptotically optimal in the minimax sense (see Arias-Castro et al. 2011 and Ingster et al. 2010). We also note that our procedure is more powerful in detecting significant variables when $|G|$ gets smaller in view of the lower bound in (13). This observation partly motivates us to consider the screening procedure in Section 3.2.

REMARK 2.4. An important byproduct of the proof of Theorem 2.4 is that the distribution of $\max_{1 \leq j \leq p} \sqrt{n} |\breve{\beta}_j - \beta_j^0| / \sqrt{\widehat{\omega}_{jj}}$ can be well approximated by $\max_{1 \leq j \leq p} |Z_j|$ with $Z = (Z_1, \ldots, Z_p) \sim^d N(0, \widetilde{\Theta})$. Therefore, we have for any $x \in \mathbb{R}$ and as $p \to +\infty$,

$$P\left(\max_{1 \leq j \leq p} n|\breve{\beta}_j - \beta_j^0|^2/\widehat{\omega}_{jj} - 2\log(p) + \log\log(p) \leq x\right) \to \exp\left\{-\frac{1}{\sqrt{\pi}} \exp\left(-\frac{x}{2}\right)\right\}. \quad (15)$$

In contrast with our method, the above alternative testing procedure has to require $p$ to diverge. In addition, the critical value obtained from the above type I extreme value distribution may not work well in practice since this weak convergence is typically slow. Instead, it is suggested to employ an "intermediate" approximation to improve the rate of convergence in the literature, e.g., Liu et al. (2008).

## 3 Applications

This section is devoted to three concrete applications of the general theoretical results developed in Section 2. Specifically, we consider (i) support recovery; (ii) testing for sparse signals; (iii) multiple testing using the step-down method.

## 3.1 Application I: Support recovery

The major goal of this section is to identify signal locations in a pre-specified set $\widetilde{G}$, i.e. support recovery. It turns out that this support recovery problem is closely related to the re-sparsifying procedure[1] applied to the de-sparsified estimator considered in van de Geer (2014) (see Lemma 2.3 therein). In comparison with Lasso, stability selection (Meinshausen and Bühlmann 2010), and the screen and clean procedure (Wasserman and Roeder 2009), simulation results in Section 5.2 illustrate that our procedure below can be more accurate in recovering signals.

Our support recovery procedure is concerned with setting a proper threshold $\tau$ in the following set

$$\widehat{\mathcal{S}}_0(\tau) = \{j \in \widetilde{G} : |\breve{\beta}_j| > \lambda_j^*(\tau)\},$$

where $\lambda_j^*(\tau) = \sqrt{\tau \widehat{\omega}_{jj} \log(p)/n}$ and $\widehat{\omega}_{jj} = \widehat{\sigma}_\epsilon^2 \widehat{\Theta}_j^T \widehat{\Sigma} \widehat{\Theta}_j$. We consider the most challenging scenario where $\widetilde{G} = [p]$ with $[p] := \{1, 2, \ldots, p\}$. In proposition 3.1, we show that the above support recovery procedure is consistent if the threshold value is set as $\tau^* = 2$, and further justify the optimality of $\tau^*$.

PROPOSITION **3.1**. *Under the assumptions in Theorem 2.3 and Assumption 2.5, we have*

$$\inf_{\beta^0 \in \Psi(2\sqrt{2})} P(\widehat{\mathcal{S}}_0(2) = \mathcal{S}_0) \to 1, \qquad (16)$$

*where* $\Psi(c_0) = \{\beta = (\beta_1, \ldots, \beta_p)^T : \min_{j \in \mathcal{S}_0} |\beta_j|/\sqrt{\omega_{jj}} > c_0 \sqrt{\log(p)/n}\}$. *Moreover, we have for any* $0 < \tau < 2$,

$$\sup_{\beta^0 \in \mathcal{S}^*(s_0)} P(\widehat{\mathcal{S}}_0(\tau) = \mathcal{S}_0) \to 0, \qquad (17)$$

*where* $\mathcal{S}^*(s_0) = \{\beta = (\beta_1, \ldots, \beta_p)^T \in \mathbb{R}^p : \sum_{j=1}^p \mathbf{I}\{\beta_j \neq 0\} = s_0\}$.

---

[1] This re-sparsifying procedure has the merits that it can improve the $l_\infty$-bounds of Lasso and has $l_q$-bounds similar to Lasso (under sparsity conditions).

A key step in the proof of Proposition 3.1 is (15) in Remark 2.4.

## 3.2 Application II: Testing for sparse signals

In this subsection, we focus on the testing problem, $H_{0,\widetilde{G}} : \beta_j^0 = 0$ for any $j \in \widetilde{G} \subseteq \{1, 2, \ldots, p\}$. To improve the efficiency of the testing procedure and reduce the computational cost in the nodewise Lasso, we propose a three-step procedure that first randomly splits the sample into two subsamples, screens out the irrelevant variables (leading to a reduced model) based on the first subsample, and then performs simultaneous testing in Section 2.2 on the reduced model based on the second subsample.

Suppose the predictors are properly centered and studentized with sample mean zero and standard deviation one. The three-step procedure is formally described as follows:

1. *Random sample splitting:* Randomly split the sample into two subsamples $\{(\widetilde{X}_i, Y_i)\}_{i \in \mathcal{D}_1}$ and $\{(\widetilde{X}_i, Y_i)\}_{i \in \mathcal{D}_2}$, where $\mathcal{D}_1 \cup \mathcal{D}_2 = \{1, 2, \ldots, n\}$, $|\mathcal{D}_1| = \lfloor c_0 n \rfloor$ and $|\mathcal{D}_2| = n - \lfloor c_0 n \rfloor$ for some $0 < c_0 < 1$.

2. *Marginal screening based on $\mathcal{D}_1$:* Let $\mathbf{X}_{\mathcal{D}_1}$ be the submatrix of $\mathbf{X}$ that contains the rows in $\mathcal{D}_1$ and let $Y_{\mathcal{D}_1} = (Y_i)_{i \in \mathcal{D}_1}$. Compute the correlation $W = (w_1, \ldots, w_p)^T = \mathbf{X}_{\mathcal{D}_1}^T Y_{\mathcal{D}_1}$ and consider the submodel $\mathcal{S}_\gamma = \{1 \leq j \leq p : |w_j| > \gamma\}$ with $\gamma$ being a positive number such that $|\mathcal{S}_\gamma| = |\mathcal{D}_2| - 1$ or $|\mathcal{S}_\gamma| = \lfloor |\mathcal{D}_2| / \log(|\mathcal{D}_2|) \rfloor$.

3. *Testing after screening based on $\mathcal{D}_2$:* Under the reduced model $\mathcal{S}_\gamma$, compute the de-sparsifying Lasso estimator $\{\breve{\beta}_j\}_{j \in \mathcal{S}_\gamma}$ and the variance estimator $\widehat{\omega}_{jj}$ based on $\{(\widetilde{X}_i, Y_i)\}_{i \in \mathcal{D}_2}$. Define $\widetilde{G}_\gamma = \widetilde{G} \cap \mathcal{S}_\gamma$. Denote by $T_{nst,\gamma} = \max_{j \in \widetilde{G}_\gamma} \sqrt{n} |\breve{\beta}_j|$ and $T_{st,\gamma} = \max_{j \in \widetilde{G}_\gamma} \sqrt{n} |\breve{\beta}_j| / \sqrt{\widehat{\omega}_{jj}}$ the non-studentized and studentized test statistics, respectively (if $\widetilde{G}_\gamma = \emptyset$, we simply set $T_{nst,\gamma} = T_{st,\gamma} = 0$ and do not reject $H_{0,\widetilde{G}}$). Let $c^*_{\widetilde{G}_\gamma}(\alpha)$ and $\bar{c}^*_{\widetilde{G}_\gamma}(\alpha)$ be the corresponding bootstrap critical values at level $\alpha$. Reject the null hypothesis if $T_{nst,\gamma} > c^*_{\widetilde{G}_\gamma}(\alpha)$ ( $T_{st,\gamma} > \bar{c}^*_{\widetilde{G}_\gamma}(\alpha)$).

In the above three step procedure, we have employed the data-splitting strategy as suggested in Wasserman and Roeder (2009) to reduce the Type I error rate due to the selection effect (the so-called selective Type I error rate). Under suitable assumptions that rule out unfaithfulness (small partial correlations), see e.g., Fan and Lv (2008), we have $\mathcal{S}_0 \subseteq \mathcal{S}_\gamma$ with an overwhelming probability, which justifies the validity of the second step. On the event that $\mathcal{S}_0 \subseteq \mathcal{S}_\gamma$, the validity and optimality of the inference procedure based on $T_{nst,\gamma}$ and $T_{st,\gamma}$ have been established in Section 2.2.

From a practical viewpoint, the three-step procedure enjoys two major advantages over the single step procedure: (1) the nodewise Lasso involves the computation of $p$ Lasso problems, which can be computationally intensive especially when $p$ is very large (e.g. $p$ could be tens of thousands in genomic studies). The three-step procedure lessens this computation burden as the nodewise Lasso is now performed under the reduced model $\mathcal{S}_\gamma$ with a much smaller size; (2) Due to the screening step, a reduced model is created which could lead to more efficient inference in some cases, see Section 5.3. This can also be seen from Theorem 2.4 that our testing procedure is more powerful when $|G|$ gets smaller.

### 3.3 Application III: Multiple testing with strong FWER control

We are interested in the following multiple testing problem:

$$H_{0,j} : \beta_j^0 \leq \widetilde{\beta}_j \text{ versus } H_{a,j} : \beta_j^0 > \widetilde{\beta}_j \text{ for all } j \in G.$$

For simplicity, we set $G = [p]$. To obtain a strong control of the FWER, we couple the bootstrap-assisted testing procedure with the step-down method proposed in Romano and Wolf (2005). Our method is a special case of a general methodology presented in Section 5 of Chernozhukov et al. (2013) by setting $\widehat{\beta}$ therein as the de-sparsifying estimator $\breve{\beta}$. As pointed out in Chernozhukov et al. (2013), our method has two important features: (i) it applies to models with an increasing dimension; (ii) it is asymptotically non-conservative as compared to the Bonferroni-Holm procedure since the correlation amongst the test statistics is taken into account. In fact, we will compare the finite sample performance of our method with that of the Bonferroni-Holm procedure in Section 5.4. We also want to point out that any procedure controlling the FWER will also control the false discovery rate (Benjamin and Hochberg 1995) when there exist some true discoveries.

Denote by $\Omega$ the space for all data generating processes, and $\omega_0$ be the true process. Each null hypothesis $H_{0,j}$ is equivalent to $\omega_0 \in \Omega_j$ for some $\Omega_j \subseteq \Omega$. For any $\eta \subseteq [p]$, denote by $\Omega^\eta = (\cap_{j \in \eta} \Omega_j) \cap (\cap_{j \notin \eta} \Omega_j^c)$ with $\Omega_j^c = \Omega \setminus \Omega_j$. The strong control of the FWER means that,

$$\sup_{\eta \subseteq [p]} \sup_{\omega_0 \in \Omega^\eta} P_{\omega_0} (\text{reject at least one hypothesis } H_{0,j}, j \in \eta) \leq \alpha + o(1). \tag{18}$$

Let $T_j = \sqrt{n}(\breve{\beta}_j - \widetilde{\beta}_j)$ and denote by $c_\eta(\alpha)$ the bootstrapped estimate for the $1 - \alpha$ quantile of $\max_{j \in \eta} T_j$. The step-down method in Romano and Wolf (2005) in controlling the FWER is described as follows. Let $\eta(1) = [p]$ at the first step. Reject all hypotheses $H_{0,j}$ such that $T_j > c_{\eta(1)}(\alpha)$. If no hypothesis is rejected, then stop. If some hypotheses are rejected, let $\eta(2)$ be the set of indices for those hypotheses not being rejected at the first step. At step $l$, let $\eta(l) \subseteq [p]$ be the subset of hypothesises that were not rejected at step $l - 1$. Reject all hypothesises $H_{0,j}, j \in \eta(l)$ satisfying that $T_j > c_{\eta(l)}(\alpha)$. If no hypothesis is rejected, then stop. Proceed in this way until the the algorithm stops. As shown in Romano and Wolf (2005), the strong control of the family-wise error holds provided that

$$c_\eta(\alpha) \leq c_{\eta'}(\alpha), \quad \text{for } \eta \subseteq \eta', \tag{19}$$

$$\sup_{\eta \subseteq [p]} \sup_{\omega_0 \in \Omega^\eta} P_{\omega_0} \left( \max_{j \in \eta} T_j > c_\eta(\alpha) \right) \leq \alpha + o(1). \tag{20}$$

Therefore, we can show that the step-down method together with the multiplier bootstrap provide strong control of the FWER by verifying (19) and (20). The arguments are similar to those in the proof of Theorem 2.2; also see Theorem 5.1 of Chernozhukov et al. (2013).

PROPOSITION **3.2**. *Under the assumptions in Theorem 2.2, the step-down procedure with the bootstrap critical value $c_\eta(\alpha)$ satisfies (18).*

## 4 Generalization

In this section, our results are extended beyond the linear models to a general framework with a convex loss function and a penalty function. For $y \in \mathcal{Y} \subseteq \mathbb{R}$ and $x \in \mathcal{X} \subseteq \mathbb{R}^p$, consider a loss function $L_\beta(y,x) = L(y, x^T\beta)$ which is strictly convex in $\beta \in \mathbb{R}^p$. The regularized estimator based on the penalty function $\rho_\lambda(\cdot)$ is defined as

$$\widehat{\beta} = \arg\min_{\beta \in \mathbb{R}^p} \left\{ \mathbb{E}_n L_\beta + \sum_{j=1}^{p} \rho_\lambda(|\beta_j|) \right\}, \tag{21}$$

where $\mathbb{E}_n g = \sum_{i=1}^{n} g(y_i, x_i)/n$ and $\{(y_i, x_i)\}_{i=1}^{n}$ is a sequence of i.i.d observations. Note that our formulation (21) slightly generalizes the framework in Section 3 of van de Geer et al. (2014) by considering a general penalty function. Again, we want to test the hypothesis $H_{0,G} : \beta_j^0 = \widetilde{\beta}_j$ for all $j \in G \subseteq [p]$ versus $H_{0,G} : \beta_j^0 \neq \widetilde{\beta}_j$ for some $j \in G$. Our test statistic is given by $\max_{j \in G} \sqrt{n}|\breve{\beta}_j - \widetilde{\beta}_j|$ with $\breve{\beta}_j$ being the de-sparsifying Lasso estimator defined below. We use analogous notation as in Section 2 but with some modifications for the current context. For example, denote $\beta^0$ as the unique minimizer of $\beta \mapsto \mathbb{E} L_\beta$.

Define

$$\dot{L}_\beta(y,x) = \frac{\partial}{\partial \beta} L_\beta(y,x) = x\dot{L}(y, x^T\beta), \quad \ddot{L}_\beta(y,x) = \frac{\partial^2}{\partial \beta \partial \beta^T} L_\beta(y,x) = xx^T \ddot{L}(y, x^T\beta),$$

where $\dot{L} = \partial L(y,a)/\partial a$ and $\ddot{L} = \partial^2 L(y,a)/\partial^2 a$. Define $\widehat{\Sigma} = \mathbb{E}_n \ddot{L}_{\widehat{\beta}}$ and let $\widehat{\Theta} := \widehat{\Theta}(\widehat{\beta})$ be a suitable approximation for the inverse of $\widehat{\Sigma}$ (see more details in Section 3.1.1 of van de Geer et al. 2014). When the penalized loss function in (21) is convex in $\beta$, the KKT condition gives

$$\mathbb{E}_n \dot{L}_{\widehat{\beta}} + \widehat{\kappa}_\lambda = 0, \tag{22}$$

where $\widehat{\kappa}_\lambda = (\widehat{\kappa}_1, \ldots, \widehat{\kappa}_p)^T$ with $\widehat{\kappa}_j = \text{sign}(\widehat{\beta}_j)\dot{\rho}_\lambda(|\widehat{\beta}_j|)$ if $\widehat{\beta}_j \neq 0$ and some $|\widehat{\kappa}_j| \leq |\dot{\rho}_\lambda(0+)|$ if $\widehat{\beta}_j = 0$. By Taylor expansion, we have

$$\mathbb{E}_n \dot{L}_{\widehat{\beta}} = \mathbb{E}_n \dot{L}_{\beta_0} + \mathbb{E}_n \ddot{L}_{\widehat{\beta}}(\widehat{\beta} - \beta^0) + \mathcal{R},$$

where $\mathcal{R}$ is the remainder term. Plugging back to the KKT condition, we obtain,

$$\mathbb{E}_n \dot{L}_{\beta_0} + \widehat{\Sigma}(\widehat{\beta} - \beta^0) + \mathcal{R} = -\widehat{\kappa}_\lambda,$$

implying that $\widehat{\beta} - \beta^0 + \Delta/\sqrt{n} + \widehat{\Theta}\widehat{\kappa}_\lambda = -\widehat{\Theta}\mathbb{E}_n \dot{L}_{\beta_0} - \widehat{\Theta}\mathcal{R}$, where $\Delta = \sqrt{n}(\widehat{\Theta}\widehat{\Sigma} - I)(\widehat{\beta} - \beta_0)$.

Following van de Geer et al. (2014), we define the de-sparsifying/de-biased estimator as

$$\breve{\beta} = \widehat{\beta} + \widehat{\Theta}\widehat{\kappa}_\lambda = \widehat{\beta} - \widehat{\Theta}\mathbb{E}_n \dot{L}_{\widehat{\beta}}.$$

Note $\breve{\beta} - \beta^0 = -\widehat{\Theta}\mathbb{E}_n \dot{L}_{\beta_0} - \widehat{\Theta}\mathcal{R} - \Delta/\sqrt{n}$. With some abuse of notation, define $\xi_{ij} = -\Theta_j^T \dot{L}_{\beta_0}(y_i, x_i) = \Theta_j^T x_i \epsilon_i$ with $\epsilon_i = -\dot{L}(y_i, x_i^T \beta_0)$, and $\widehat{\xi}_{ij} = -\widehat{\Theta}_j^T \dot{L}_{\beta_0}(y_i, x_i)$. To conduct inference, we employ the multiplier bootstrap in the following way. Generate a sequence of i.i.d standard normal random variables $\{e_i\}$ and define the bootstrap statistic,

$$\widetilde{W}_G^* = \max_{j \in G} \left| \sum_{i=1}^n \widehat{\Theta}_j^T x_i \dot{L}(y_i, x_i^T \widetilde{\beta}) e_i / \sqrt{n} \right|,$$

where $\widetilde{\beta}$ is a suitable estimator for $\beta^0$. The bootstrap critical value is given by $\widetilde{c}_G^*(\alpha) = \inf\{t \in \mathbb{R} : P(\widetilde{W}_G^* \leq t | \{(y_i, x_i)\}_{i=1}^n) \geq 1 - \alpha\}$.

We first impose Assumptions 4.1–4.4 that are similar to Assumptions (C1)-(C5) in van de Geer et al. (2014). Denote by $\mathbf{X}$ the design matrix with the $i$th row equal to $x_i^T$.

ASSUMPTION 4.1. *The derivatives $\dot{L}(y, a)$ and $\ddot{L}(y, a)$ exist for all $y$ and $a$. For some $\delta$-neighborhood,*

$$\max_{a_0 \in \{x_i^T \beta^0 : x_i \in \mathcal{X}\}} \sup_{|a - a_0| \vee |a' - a_0| \leq \delta} \sup_{y \in \mathcal{Y}} \frac{|\ddot{L}(y, a) - \ddot{L}(y, a')|}{|a - a'|} \leq 1.$$

ASSUMPTION 4.2. *Assume $\max_{i,j}|x_{ij}| \leq K_n$.*

ASSUMPTION 4.3. *Assume that $||\widehat{\beta} - \beta^0||_1 = O_P(s_0\lambda)$ and $||\mathbf{X}(\widehat{\beta} - \beta^0)||_2^2/n = O_P(s_0\lambda^2)$.*

ASSUMPTION 4.4. *Assume that $||\mathbb{E}_n \ddot{L}_{\widehat{\beta}} \widehat{\Theta}_j - 1_j||_\infty = O_P(\lambda_*)$ for some $\lambda_* > 0$, and $||\mathbf{X}\widehat{\Theta}_j||_\infty = O_P(K_n)$ uniformly for $j$. Here, $1_j$ denotes the vector with the $j$-th element one and others zero.*

Moreover, we make the following additional assumptions. In particular, Assumptions 4.5-4.6 are parallel to Assumptions 2.3-2.4, while Assumption 4.7 is a technical one motivated by the results in Theorem 3.2 of van de Geer et al. (2014).

ASSUMPTION 4.5. *Suppose $\max_{1 \leq j \leq p} s_j K_n^2 (\log(pn))^7/n \leq C_2 n^{-c_2}$ for some constants $c_2, C_2 > 0$. Suppose $\{\epsilon_i\}$ are i.i.d sub-exponential that is $\mathbb{E}\exp(|\epsilon_i|/C') \leq 1$ for some large enough constants $C' > 0$, and $c' < \sigma_\epsilon^2 < c$ for $c, c' > 0$.*

Let $\widetilde{\Gamma} = \max_{1 \leq j,k \leq p} |\widehat{\Theta}_j^T \widehat{\Sigma}_{\widetilde{\beta}} \widehat{\Theta}_k - \Theta_j^T \Sigma_{\beta_0} \Theta_k|$, where $\widehat{\Sigma}_{\widetilde{\beta}} = \mathbf{X}^T \mathbf{W}_{\widetilde{\beta}} \mathbf{X}/n$ with $\mathbf{W}_{\widetilde{\beta}} = \mathrm{diag}(\dot{L}^2(y_1, x_1^T \widetilde{\beta}), \ldots, \dot{L}^2(y_n, x_n^T \widetilde{\beta}))$, and $\Sigma_{\beta_0} = \mathrm{cov}(x_i \dot{L}(y_i, x_i^T \beta_0))$.

ASSUMPTION 4.6. *There exists a sequence of positive numbers $\alpha_n \to +\infty$ such that $\alpha_n/p = o(1)$ and $P(\alpha_n (\log p)^2 \widetilde{\Gamma} > 1) \to 0$.*

ASSUMPTION 4.7. *Assume that uniformly for $j$, it holds that*

$$||\widehat{\Theta}_j(\widehat{\beta}) - \Theta_j(\beta^0)||_1 = O_P(K_n s_j \sqrt{\log(p)/n}) + O_P\left(K_n^2 s_0((\lambda^2/\sqrt{\log(p)/n}) \vee \lambda)\right).$$

We are now in position to present the main result in this section which justifies the use of the bootstrap procedure.

THEOREM 4.1. *Suppose Assumptions 4.1-4.7 hold. Suppose*

$$\sqrt{n \log(p)} K_n s_0 \lambda^2 = o(1), \quad \sqrt{n \log(p)} \lambda \lambda_* s_0 = o(1), \quad \max_j s_j K_n (\log(p))^{3/2}/\sqrt{n} = o(1),$$

$$\sqrt{\log(p)} K_n^2 s_0 \left(\lambda^2 \sqrt{n} \vee \lambda \sqrt{\log(p)}\right) = o(1), \quad K_n^2 \log(p)(\log(n))^2/n = o(1).$$

Then we have for any $G \subseteq \{1, 2, \ldots, p\}$,

$$\sup_{\alpha \in (0,1)} \left| P\left(\max_{j \in G} \sqrt{n}|\breve{\beta}_j - \beta_j^0| > \widetilde{c}_G^*(\alpha)\right) - \alpha \right| = o(1). \tag{23}$$

REMARK 4.1. Let $\widetilde{W}_G^{**} = \max_{j \in G} \left| \sum_{i=1}^n \widehat{\Theta}_j^T x_i \dot{L}(y_i, x_i^T \widetilde{\beta}) e_i / \sqrt{n \widehat{w}_{jj}} \right|$, where $\omega_{jj} = \widehat{\Theta}_j^T \widehat{\Sigma}_{\widetilde{\beta}} \widehat{\Theta}_j$ and $\{e_i\}$ is a sequence of i.i.d standard normal random variables and $\widehat{w}_{jj} = \widehat{\Theta}_j^T \widehat{\Sigma}_{\widetilde{\beta}} \widehat{\Theta}_j$. Define the bootstrap critical value $\widetilde{c}_G^{**}(\alpha) = \inf\{t \in \mathbb{R} : P(\widetilde{W}_G^{**} \leq t | \{(y_i, x_i)\}_{i=1}^n) \geq 1 - \alpha\}$. Following the arguments in the proofs of Theorem 2.3 and Theorem 4.1, we expect a similar result for the studentized test statistic $\max_{j \in G} \sqrt{n} |(\breve{\beta}_j - \beta_j^0) / \sqrt{\widehat{\omega}_{jj}}|$. The technical details are omitted here.

# 5 Numerical Results

In this section, we conduct some simulation studies to evaluate the finite sample performance of the methods proposed in Section 3. All the results are obtained based on sample size $n = 100$, and $1,000$ Monte Carlo replications.

To obtain the main Lasso estimator, we implemented the scaled Lasso with the tuning parameter $\lambda_0 = \sqrt{2} \tilde{L}_n(k_0/p)$ with $\tilde{L}_n(t) = n^{-1/2} \Phi^{-1}(1-t)$, where $\Phi$ is the cumulative distribution function for $N(0,1)$, and $k_0$ is the solution to $k = \tilde{L}_1^4(k/p) + 2\tilde{L}_1^2(k/p)$ (see Sun and Zhang 2013). We estimate the noise level $\sigma^2$ using the modified variance estimator defined in (24) below. The tuning parameters $\lambda_j$s in the nodewise Lasso are chosen via 10-fold cross-validation among all nodewise regressions throughout the simulation studies.

## 5.1 Simultaneous confidence intervals

We consider the linear models where the rows of **X** are fixed i.i.d realizations from $N_p(0, \Sigma)$ with $\Sigma = (\Sigma_{i,j})_{i,j=1}^p$ under two scenarios: (i) Toeplitz: $\Sigma_{i,j} = 0.9^{|i-j|}$; (ii) Exchangeable/Compound symmetric: $\Sigma_{i,i} = 1$ and $\Sigma_{i,j} = 0.8$ for $i \neq j$. The active set is $\{1, 2, \ldots, s_0\}$,

where $s_0 = |\mathcal{S}_0| = 3$ or 15. The coefficients of the linear regression models are generated according to Unif$[0, 2]$ (uniform distribution on $[0, 2]$), and the errors are generated from the studentized $t(4)$ distribution, i.e., $t(4)/\sqrt{2}$, or centeralized and the studentized Gamma(4,1) distribution i.e., (Gamma$(4, 1) - 4)/2$.

In our simulations, we found that the scaled Lasso tends to underestimate the noise level which could lead to undercoverage. To overcome this problem, we suggest the following modified variance estimator (see a similar estimator in Reid et al. 2014),

$$\widehat{\sigma}^2 = \frac{1}{n - ||\widehat{\beta}_{\text{sc}}||_0}||Y - \mathbf{X}\widehat{\beta}_{\text{sc}}||_2^2, \tag{24}$$

where $\widehat{\beta}_{\text{sc}}$ denotes the scaled Lasso estimator with the tuning parameter $\lambda_0$. Figure ?? in the supplement provides boxplots of $\widehat{\sigma}/\sigma$ for the variance estimator delivered by the scaled Lasso (denoted by "SLasso") and for the modified variance estimator in (24) (denoted by "SLasso*"). Clearly, the modified variance estimator corrects the noise underestimation issue, and thus is more preferable.

In Tables 1-2, we present the coverage probabilities and interval widths for the simultaneous confidence intervals in three different cases: $G = \mathcal{S}_0, \mathcal{S}_0^c$ or $[p]$ (a similar setup was considered in van de Geer et al. 2014). For each simulation run, we record whether the simultaneous confidence interval contains $\beta_j^0$ for $j \in G$ and the corresponding interval width. The coverage probabilities and interval widths are then calculated by averaging over 1,000 simulation runs. It is not surprising that the coverage probability is affected by the dimension $p$, the tuning parameters $\lambda_j$ (in the nodewise Lasso), the cardinality $s_0$, and the covariance matrix $\Sigma$. Overall, the non-studentized method provides satisfactory coverage probability. However, the method based on the extreme value distribution approximation is invalid when the dimension of the components of interest is low; see Remark 2.4. To avoid sending a misleading message, we choose not to provide the results based on the extreme value distribution approximation when $G = \mathcal{S}_0$.

More specifically, we observe from Tables 1-2 that: (i) the coverage is in general more accurate for $\mathcal{S}_0^c$ (as compared to $\mathcal{S}_0$), even though $|\mathcal{S}_0^c| \gg |\mathcal{S}_0| = s_0$. Similar finding is found in van de Geer et al. (2014) where the coverage for a single coefficient in $\mathcal{S}_0^c$ is more accurate; (ii) the non-studentized test statistic tends to provide better coverage (but with larger width) as compared to its studentized version when $s_0$ is large; (iii) when $s_0 = 15$, it becomes difficult to provide accurate coverage for all the active variables among the candidates. In this case, the coverage for the active set can be significantly lower than the nominal level. Additional numerical results in the supplement (see Figure ??) indicate that the undercoverage in this case is closely related with the remainder term whose maximum norm generally increases with $s_0$; (iv) compared to the confidence intervals for individual coefficients in van de Geer et al. (2014), the interval widths for simultaneous intervals are wider, which reflects the price we pay for multiplicity adjustment.

## 5.2 Support recovery

Below we compare the finite sample performance of the support recovery procedure described in Section 3.1 with those of Lasso, stability selection (Meinshausen and Bühlmann, 2010), and the screen and clean procedure in Wasserman and Roeder (2009). Consider the simulation setup in Section 5.1, where the coefficients are now generated from Unif $[2, 4]$, and the support set $\mathcal{S}_0 = \{u_1, \ldots, u_{s_0}\}$ with $u_1, \ldots, u_{s_0}$ being a realization of $s_0$ i.i.d draws without replacement from $\{1, \ldots, p\}$. To assess the performance, we consider the following similarity measure

$$d(\widehat{\mathcal{S}}_0, \mathcal{S}_0) = \frac{|\widehat{\mathcal{S}}_0 \cap \mathcal{S}_0|}{\sqrt{|\widehat{\mathcal{S}}_0| \cdot |\mathcal{S}_0|}}.$$

In the implementation of stability selection, we choose the threshold for selection frequency to be 0.6 and the upper bound for the expected number of false positives to be 2.5; see Meinshausen and Bühlmann (2010). For the screen and clean procedure, we randomly split the data into two groups, conduct screening on the first half of the data and cleaning on the second half. This two splits procedure was advocated in the simulations of Wasserman and Roeder (2009). Table 3 summarizes the mean and standard deviation of $d(\widehat{\mathcal{S}}_0, \mathcal{S}_0)$ as well as the numbers of false positives (FP) and false negatives (FN) based on 1,000 simulation runs. When $s_0 = 3$, the proposed support recovery procedure clearly outperforms Lasso, and it is comparable to stability selection and the screen and clean procedure. When $s_0 = 15$, the recovery procedure in general outperforms all other three competitors. We note that when $s_0 = 15$ and $\Sigma$ is exchangeable, the stability selection gives a high number of FN, and thus results in a low value of $d(\widehat{\mathcal{S}}_0, \mathcal{S}_0)$. This is not surprising given that stability selection is mainly designed for conservatively controlling the expected number of FP. The upper bound 2.5 set above might be too small for $s_0 = 15$. The screen and clean procedure generally performs well for $p = 120$, but its performance deteriorates for $p = 500$ and $s_0 = 15$. Overall, the proposed method performs quite well as compared to some existing alternatives.

## 5.3 Testing for sparse signals

This subsection is devoted to empirically examine the three-step testing procedure proposed in Section 3.2. We consider the following two scenarios: (i) $\Sigma$ is Toeplitz with $\Sigma_{i,j} = 0.9^{|i-j|}$, and $\beta_j^0 = \sqrt{10 \log(p)/n}$ for $1 \leq j \leq s_0$ and zero otherwise; (ii) $\Sigma$ is exchangeable, i.e., $\Sigma_{i,i} = 1$ and $\Sigma_{i,j} = 0.8$ for $i \neq j$, and $\beta_j^0 = 10\sqrt{\log(p)/n}$ for $1 \leq j \leq s_0$ and zero otherwise.

To implement the three-step testing procedure, we choose the splitting proportion $c_0 = 1/5$ in case (i) and $c_0 = 1/3$ in case (ii). For the marginal screening step, we pick $\gamma$ such that $|\mathcal{S}_\gamma| = |\mathcal{D}_2| - 1$. As pointed out in Fan and Lv (2008), the marginal screening may not perform very well in the case of strong pairwise correlation, i.e., when $\Sigma$ is exchangeable. To overcome this problem, we propose the following remedy inspired by the iterative sure independence screening in Fan and Lv (2008). We first select a subset $\mathcal{B}_1$ of $k_1$ variables using Lasso based on the subsample $\mathcal{D}_1$. Let $\widehat{r}_i = Y_i - \widetilde{X}_i^T \widehat{\beta}$ be the corresponding residuals with $\widehat{\beta}$ being the Lasso estimator, and $\widetilde{X}_i^c = (X_{ij})_{j \notin \mathcal{B}_1}$. We then treat those residuals as the new responses and apply the marginal screening as described in Step 2 to $\{(\widetilde{X}_i^c, \widehat{r}_i)\}_{i \in \mathcal{D}_1}$ to select a subset $\mathcal{B}_2$ of $|\mathcal{D}_2| - 1 - k_1$ variables. Finally, we let $\mathcal{B} = \mathcal{B}_1 \cup \mathcal{B}_2$, which contains $|\mathcal{D}_2| - 1$ variables. In case (ii) with $s_0 = 3$, this remedy selects all the three relevant variables with probability 0.98 which is much higher than 0.59 delivered by the marginal screening.

The numerical results are reported in Table 4, which compares the performance between the three-step procedure and the one-step procedure without sample splitting and marginal screening. Some remarks are in order regarding the simulation results: (1) for the Toeplitz covariance structure, the three-step procedure has reasonable size for $t$ errors, and its size is slightly upward distorted for Gamma errors. In contrast, the one-step procedure has downward size distortion; for exchangeable covariance structure, both procedures show upward size distortions; (2) the three-step procedure generates higher power for the Toeplitz covariance structure. We note that the empirical powers of both procedures are close to 1 in case (ii) due to the higher signal-to-noise ratio. Overall, the three-step procedure has better power property in case (i).

## 5.4 Testing with FWER control

In this section, we compare the finite sample performance of the step-down method in Section 3.3 with the Bonferroni-Holm procedure whose finite sample performance has been studied in van de Geer et al. (2014). To this end, we consider multiple two sided testing of hypothesis $H_{0,j} : \beta_j^0 = 0$ among all $j = 1, 2, \ldots, p$. The data is again generated from the linear model considered in Section 5.1. For $\Sigma$, we focus on the following two cases: (i) Toeplitz: $\Sigma_{i,j} = 0.9^{|i-j|}$; (ii) Block diagonal: $\Sigma_{i,i} = 1$, $\Sigma_{i,j} = 0.9$ for $5(k-1) + 1 \leq i \neq j \leq 5k$ with $k = 1, 2, \ldots, \lfloor p/5 \rfloor$, and $\Sigma_{i,j} = 0$ otherwise. We employ both the step-down method based on the studentized/non-studentized test statistic, and the Bonferroni-Holm procedure (based on the studentized test statistic) to control the FWER. Table 5 reports both the FWER and the average power, which is defined as $\sum_{j \in \mathcal{S}_0} \mathbf{I}\{H_{0,j} \text{ is rejected}\}/s_0$ based on $1,000$ simulated data sets. As seen from the table, the two procedures provide similar control on the FWER.

And the step-down method delivers slightly higher average power across all cases considered here.

## 6 Conclusion

In this paper, we have introduced a bootstrap-assisted procedure to conduct simultaneous inference for high-dimensional components of a large parameter vector in sparse linear models. Our procedure is proved to achieve the pre-specified significance level asymptotically and to enjoy certain optimality in terms of its power. Our general theory has been successfully applied to three concrete examples, namely support recovery, testing for sparse signals, and multiple testing using the step-down method. Below we point out a few future research directions. The first direction is the automatic and efficient selection of the tuning parameters (e.g. $\lambda_j$s in the nodewise Lasso) in the context of hypothesis testing and confidence interval construction (see e.g. Appendix A.1 of Dezeure et al. 2014). The second direction is to adapt the proposed method to conduct inference on high dimensional concentration matrix (see e.g. Jankova and van de Geer 2014). Finally, it is also interesting to extend our results to more complicated models e.g. Cox model, after some technical modifications.

## Acknowledgments

The authors would like to thank the Editor, Associate Editor and reviewers for their constructive comments and helpful suggestions, which substantially improved the paper.

Table 1: Coverage probabilities and interval widths for the simultaneous confidence intervals based on the non-studentized ("NST") and studentized ("ST") test statistics, where $\mathcal{S}_0 = \{1, 2, 3\}$ and $p = 120, 500$. The row "EX" corresponds to the coverage based on the studentized test statistic with the Type I extreme distribution approximation. "Cov" and "Len" denote the coverage probability and interval width respectively. Case (i) ((ii)) corresponds to Toeplitz matrix (exchangeable matrix).

| | | | $\mathcal{S}_0$,(i) | | $\mathcal{S}_0^c$,(i) | | $[p]$,(i) | | $\mathcal{S}_0$,(ii) | | $\mathcal{S}_0^c$,(ii) | | $[p]$,(ii) | |
| --- | --- | --- | --- | --- | --- | --- | --- | --- | --- | --- | --- | --- | --- | --- |
| | | | 95% | 99% | 95% | 99% | 95% | 99% | 95% | 99% | 95% | 99% | 95% | 99% |
| $p = 120$ | | | | | | | | | | | | | | |
| $t(4)/\sqrt{2}$ | $\text{NST}_{\text{cv}}$ | Cov | 0.82 | 0.94 | 0.97 | 0.99 | 0.95 | 0.99 | 0.91 | 0.97 | 0.93 | 0.98 | 0.93 | 0.98 |
| | | Len | 0.99 | 1.22 | 1.49 | 1.67 | 1.50 | 1.67 | 0.97 | 1.18 | 1.49 | 1.67 | 1.49 | 1.67 |
| | $\text{ST}_{\text{cv}}$ | Cov | 0.82 | 0.93 | 0.97 | 0.99 | 0.96 | 0.99 | 0.92 | 0.97 | 0.92 | 0.98 | 0.92 | 0.98 |
| | | Len | 0.97 | 1.19 | 1.48 | 1.64 | 1.48 | 1.64 | 0.96 | 1.18 | 1.46 | 1.62 | 1.46 | 1.63 |
| | $\text{EX}_{\text{cv}}$ | Cov | NA | NA | 0.98 | 1.00 | 0.97 | 0.99 | NA | NA | 0.94 | 0.98 | 0.93 | 0.98 |
| | | Len | NA | NA | 1.51 | 1.69 | 1.51 | 1.69 | NA | NA | 1.49 | 1.66 | 1.49 | 1.67 |
| Gamma | $\text{NST}_{\text{cv}}$ | Cov | 0.84 | 0.93 | 0.96 | 0.99 | 0.94 | 0.98 | 0.91 | 0.97 | 0.93 | 0.98 | 0.93 | 0.98 |
| | | Len | 0.99 | 1.22 | 1.50 | 1.67 | 1.50 | 1.67 | 0.97 | 1.18 | 1.50 | 1.68 | 1.50 | 1.68 |
| | $\text{ST}_{\text{cv}}$ | Cov | 0.82 | 0.92 | 0.97 | 0.99 | 0.95 | 0.98 | 0.90 | 0.97 | 0.92 | 0.98 | 0.92 | 0.98 |
| | | Len | 0.97 | 1.19 | 1.48 | 1.65 | 1.48 | 1.65 | 0.97 | 1.18 | 1.46 | 1.63 | 1.46 | 1.63 |
| | $\text{EX}_{\text{cv}}$ | Cov | NA | NA | 0.97 | 0.99 | 0.96 | 0.99 | NA | NA | 0.93 | 0.99 | 0.93 | 0.99 |
| | | Len | NA | NA | 1.51 | 1.69 | 1.51 | 1.69 | NA | NA | 1.49 | 1.67 | 1.49 | 1.67 |
| $p = 500$ | | | | | | | | | | | | | | |
| $t(4)/\sqrt{2}$ | $\text{NST}_{\text{cv}}$ | Cov | 0.76 | 0.90 | 0.96 | 0.99 | 0.94 | 0.98 | 0.92 | 0.98 | 0.92 | 0.97 | 0.92 | 0.97 |
| | | Len | 0.89 | 1.09 | 1.47 | 1.62 | 1.47 | 1.62 | 0.97 | 1.19 | 1.65 | 1.82 | 1.65 | 1.82 |
| | $\text{ST}_{\text{cv}}$ | Cov | 0.77 | 0.90 | 0.97 | 0.99 | 0.95 | 0.98 | 0.92 | 0.97 | 0.92 | 0.97 | 0.91 | 0.97 |
| | | Len | 0.88 | 1.08 | 1.46 | 1.60 | 1.46 | 1.60 | 0.97 | 1.18 | 1.62 | 1.77 | 1.62 | 1.77 |
| | $\text{EX}_{\text{cv}}$ | Cov | NA | NA | 0.98 | 0.99 | 0.96 | 0.98 | NA | NA | 0.92 | 0.97 | 0.92 | 0.97 |
| | | Len | NA | NA | 1.48 | 1.63 | 1.48 | 1.63 | NA | NA | 1.64 | 1.81 | 1.64 | 1.81 |
| Gamma | $\text{NST}_{\text{cv}}$ | Cov | 0.77 | 0.90 | 0.98 | 0.99 | 0.96 | 0.98 | 0.91 | 0.97 | 0.95 | 0.98 | 0.94 | 0.98 |
| | | Len | 0.90 | 1.10 | 1.49 | 1.63 | 1.49 | 1.64 | 0.98 | 1.19 | 1.66 | 1.83 | 1.66 | 1.83 |
| | $\text{ST}_{\text{cv}}$ | Cov | 0.77 | 0.90 | 0.98 | 0.99 | 0.96 | 0.98 | 0.90 | 0.97 | 0.93 | 0.97 | 0.93 | 0.97 |
| | | Len | 0.89 | 1.09 | 1.47 | 1.61 | 1.47 | 1.61 | 0.97 | 1.19 | 1.63 | 1.78 | 1.63 | 1.78 |
| | $\text{EX}_{\text{cv}}$ | Cov | NA | NA | 0.98 | 1.00 | 0.96 | 0.99 | NA | NA | 0.94 | 0.98 | 0.94 | 0.98 |
| | | Len | NA | NA | 1.49 | 1.64 | 1.49 | 1.64 | NA | NA | 1.65 | 1.82 | 1.65 | 1.82 |

Note: The tuning parameters $\lambda_j$s in the nodewise Lasso are chosen to be the same via 10-fold cross-validation among all nodewise regressions for $\text{NST}_{\text{cv}}$, $\text{ST}_{\text{cv}}$, and $\text{EX}_{\text{cv}}$. $t(4)/\sqrt{2}$ and Gamma denote the studentized $t(4)$ distribution and the centralized and studentized Gamma(4,1) distribution respectively. The coverage probabilities and interval widths are computed based on 1,000 simulation runs. For the studentized test, we report the average interval widths over different components.

Table 2: Coverage probabilities and interval widths for the simultaneous confidence intervals based on the non-studentized ("NST") and studentized ("ST") test statistics, where $\mathcal{S}_0 = \{1, 2, 3, \ldots, 15\}$ and $p = 120, 500$. The row "EX" corresponds to the coverage based on the studentized test statistic with the Type I extreme distribution approximation. "Cov" and "Len" denote the coverage probability and interval width respectively. Case (i) ((ii)) corresponds to Toeplitz matrix (exchangeable matrix).

| | | | $\mathcal{S}_0$,(i) | | $\mathcal{S}_0^c$,(i) | | $[p]$,(i) | | $\mathcal{S}_0$,(ii) | | $\mathcal{S}_0^c$,(ii) | | $[p]$,(ii) | |
|---|---|---|---|---|---|---|---|---|---|---|---|---|---|---|
| | | | 95% | 99% | 95% | 99% | 95% | 99% | 95% | 99% | 95% | 99% | 95% | 99% |
| $p = 120$ | | | | | | | | | | | | | | |
| $t(4)/\sqrt{2}$ | NST$_{cv}$ | Cov | 0.68 | 0.87 | 0.99 | 1.00 | 0.87 | 0.95 | 0.73 | 0.88 | 0.92 | 0.98 | 0.89 | 0.96 |
| | | Len | 1.44 | 1.73 | 1.68 | 1.92 | 1.72 | 1.97 | 1.27 | 1.48 | 1.67 | 1.90 | 1.67 | 1.91 |
| | ST$_{cv}$ | Cov | 0.50 | 0.70 | 0.99 | 1.00 | 0.75 | 0.85 | 0.68 | 0.84 | 0.83 | 0.92 | 0.74 | 0.88 |
| | | Len | 1.30 | 1.50 | 1.55 | 1.73 | 1.56 | 1.74 | 1.23 | 1.42 | 1.53 | 1.70 | 1.53 | 1.71 |
| | EX$_{cv}$ | Cov | NA | NA | 0.99 | 1.00 | 0.77 | 0.88 | NA | NA | 0.85 | 0.94 | 0.77 | 0.92 |
| | | Len | NA | NA | 1.58 | 1.77 | 1.59 | 1.78 | NA | NA | 1.55 | 1.74 | 1.56 | 1.75 |
| Gamma | NST$_{cv}$ | Cov | 0.69 | 0.88 | 0.99 | 1.00 | 0.86 | 0.95 | 0.77 | 0.90 | 0.93 | 0.98 | 0.90 | 0.97 |
| | | Len | 1.44 | 1.74 | 1.69 | 1.93 | 1.72 | 1.97 | 1.27 | 1.49 | 1.68 | 1.91 | 1.68 | 1.92 |
| | ST$_{cv}$ | Cov | 0.52 | 0.70 | 0.99 | 1.00 | 0.74 | 0.86 | 0.71 | 0.87 | 0.85 | 0.94 | 0.78 | 0.91 |
| | | Len | 1.30 | 1.51 | 1.55 | 1.73 | 1.57 | 1.75 | 1.24 | 1.43 | 1.53 | 1.71 | 1.54 | 1.71 |
| | EX$_{cv}$ | Cov | NA | NA | 0.99 | 1.00 | 0.77 | 0.88 | NA | NA | 0.87 | 0.96 | 0.81 | 0.94 |
| | | Len | NA | NA | 1.59 | 1.78 | 1.60 | 1.79 | NA | NA | 1.56 | 1.75 | 1.57 | 1.76 |
| $p = 500$ | | | | | | | | | | | | | | |
| $t(4)/\sqrt{2}$ | NST$_{cv}$ | Cov | 0.56 | 0.82 | 0.99 | 1.00 | 0.86 | 0.95 | 0.36 | 0.57 | 0.93 | 0.98 | 0.81 | 0.94 |
| | | Len | 1.35 | 1.65 | 1.74 | 1.95 | 1.75 | 1.97 | 1.32 | 1.55 | 2.01 | 2.38 | 2.01 | 2.38 |
| | ST$_{cv}$ | Cov | 0.34 | 0.57 | 0.99 | 1.00 | 0.76 | 0.86 | 0.35 | 0.55 | 0.61 | 0.77 | 0.49 | 0.66 |
| | | Len | 1.22 | 1.40 | 1.58 | 1.73 | 1.59 | 1.74 | 1.28 | 1.48 | 1.71 | 1.87 | 1.71 | 1.87 |
| | EX$_{cv}$ | Cov | NA | NA | 0.99 | 1.00 | 0.78 | 0.87 | NA | NA | 0.65 | 0.81 | 0.53 | 0.70 |
| | | Len | NA | NA | 1.61 | 1.77 | 1.61 | 1.77 | NA | NA | 1.73 | 1.90 | 1.73 | 1.91 |
| Gamma | NST$_{cv}$ | Cov | 0.52 | 0.80 | 1.00 | 1.00 | 0.86 | 0.95 | 0.38 | 0.57 | 0.94 | 0.99 | 0.83 | 0.96 |
| | | Len | 1.36 | 1.66 | 1.75 | 1.96 | 1.76 | 1.98 | 1.34 | 1.57 | 2.04 | 2.42 | 2.04 | 2.42 |
| | ST$_{cv}$ | Cov | 0.32 | 0.52 | 0.99 | 1.00 | 0.74 | 0.85 | 0.34 | 0.53 | 0.60 | 0.77 | 0.47 | 0.66 |
| | | Len | 1.22 | 1.41 | 1.59 | 1.75 | 1.60 | 1.75 | 1.29 | 1.50 | 1.73 | 1.90 | 1.73 | 1.90 |
| | EX$_{cv}$ | Cov | NA | NA | 0.99 | 1.00 | 0.76 | 0.87 | NA | NA | 0.64 | 0.80 | 0.50 | 0.69 |
| | | Len | NA | NA | 1.61 | 1.78 | 1.62 | 1.78 | NA | NA | 1.75 | 1.93 | 1.76 | 1.93 |

Note: The tuning parameters $\lambda_j$s in the nodewise Lasso are chosen via 10-fold cross-validation among all nodewise regressions for NST$_{cv}$, ST$_{cv}$, and EX$_{cv}$. $t(4)/\sqrt{2}$ and Gamma denote the studentized $t(4)$ distribution and the centralized and studentized Gamma(4,1) distribution respectively. The coverage probabilities and interval widths are computed based on 1,000 simulation runs. For the studentized test, we report the average interval widths over different components.

Table 3: The mean and standard deviation (SD) of $d(\widehat{\mathcal{S}}_0, \mathcal{S}_0)$, and the numbers of false positives (FP) and false negatives (FN). Case (i) ((ii)) corresponds to Toeplitz matrix (exchangeable matrix).

|  |  | $s_0 = 3$ | | | | $s_0 = 15$ | | | |
| --- | --- | --- | --- | --- | --- | --- | --- | --- | --- |
|  |  | Mean | SD | FP | FN | Mean | SD | FP | FN |
| $p = 120$, (i) | | | | | | | | | |
| $t(4)/\sqrt{2}$ | SupRec | 0.98 | 0.05 | 0.16 | 0.00 | 0.98 | 0.03 | 0.68 | 0.02 |
|  | Stability | 0.93 | 0.08 | 0.52 | 0.00 | 0.35 | 0.07 | 0.63 | 12.64 |
|  | Lasso$_{sc}$ | 0.68 | 0.10 | 3.94 | 0.00 | 0.72 | 0.04 | 13.87 | 0.00 |
|  | S&C | 0.99 | 0.04 | 0.04 | 0.01 | 0.87 | 0.15 | 0.18 | 3.14 |
| Gamma | SupRec | 0.97 | 0.06 | 0.20 | 0.00 | 0.98 | 0.03 | 0.71 | 0.00 |
|  | Stability | 0.93 | 0.08 | 0.58 | 0.00 | 0.36 | 0.07 | 0.61 | 12.63 |
|  | Lasso$_{sc}$ | 0.68 | 0.11 | 3.98 | 0.00 | 0.72 | 0.04 | 13.94 | 0.00 |
|  | S&C | 0.99 | 0.03 | 0.03 | 0.01 | 0.87 | 0.15 | 0.14 | 3.30 |
| $p = 120$, (ii) | | | | | | | | | |
| $t(4)/\sqrt{2}$ | SupRec | 0.97 | 0.06 | 0.24 | 0.00 | 0.98 | 0.02 | 0.55 | 0.01 |
|  | Stability | 0.97 | 0.06 | 0.26 | 0.00 | 0.00 | 0.00 | 0.13 | 15.00 |
|  | Lasso$_{sc}$ | 0.56 | 0.09 | 7.09 | 0.00 | 0.65 | 0.04 | 20.24 | 0.00 |
|  | S&C | 0.99 | 0.04 | 0.04 | 0.01 | 0.80 | 0.24 | 0.16 | 4.41 |
| Gamma | SupRec | 0.97 | 0.06 | 0.23 | 0.00 | 0.98 | 0.02 | 0.60 | 0.00 |
|  | Stability | 0.96 | 0.06 | 0.27 | 0.00 | 0.00 | 0.00 | 0.13 | 15.00 |
|  | Lasso$_{sc}$ | 0.56 | 0.09 | 7.07 | 0.00 | 0.65 | 0.04 | 20.20 | 0.00 |
|  | S&C | 0.99 | 0.03 | 0.04 | 0.00 | 0.81 | 0.22 | 0.24 | 4.40 |
| $p = 500$, (i) | | | | | | | | | |
| $t(4)/\sqrt{2}$ | SupRec | 0.97 | 0.06 | 0.20 | 0.00 | 0.53 | 0.04 | 0.14 | 10.66 |
|  | Stability | 0.84 | 0.10 | 1.32 | 0.00 | 0.36 | 0.03 | 1.74 | 11.97 |
|  | Lasso$_{sc}$ | 0.62 | 0.09 | 5.24 | 0.00 | 0.38 | 0.03 | 19.25 | 7.37 |
|  | S&C | 0.96 | 0.13 | 0.36 | 0.12 | 0.12 | 0.17 | 18.26 | 14.04 |
| Gamma | SupRec | 0.97 | 0.06 | 0.18 | 0.00 | 0.53 | 0.04 | 0.15 | 10.68 |
|  | Stability | 0.86 | 0.09 | 1.24 | 0.00 | 0.36 | 0.03 | 1.75 | 11.93 |
|  | Lasso$_{sc}$ | 0.62 | 0.09 | 5.38 | 0.00 | 0.38 | 0.03 | 19.15 | 7.36 |
|  | S&C | 0.96 | 0.14 | 0.55 | 0.14 | 0.56 | 0.08 | 18.22 | 14.00 |
| $p = 500$, (ii) | | | | | | | | | |
| $t(4)/\sqrt{2}$ | SupRec | 0.97 | 0.06 | 0.25 | 0.00 | 0.96 | 0.04 | 1.38 | 0.10 |
|  | Stability | 0.96 | 0.07 | 0.34 | 0.00 | 0.08 | 0.10 | 0.86 | 14.59 |
|  | Lasso$_{sc}$ | 0.41 | 0.05 | 15.16 | 0.00 | 0.50 | 0.02 | 45.64 | 0.00 |
|  | S&C | 0.98 | 0.08 | 0.05 | 0.05 | 0.04 | 0.11 | 28.51 | 14.61 |
| Gamma | SupRec | 0.97 | 0.07 | 0.27 | 0.00 | 0.96 | 0.04 | 1.43 | 0.04 |
|  | Stability | 0.95 | 0.07 | 0.36 | 0.00 | 0.08 | 0.10 | 0.86 | 14.58 |
|  | Lasso$_{sc}$ | 0.41 | 0.05 | 15.18 | 0.00 | 0.50 | 0.02 | 45.41 | 0.00 |
|  | S&C | 0.98 | 0.06 | 0.05 | 0.05 | 0.05 | 0.11 | 27.54 | 14.62 |

Note: The tuning parameters $\lambda_j$s in the nodewise Lasso are chosen to be the same via 10-fold cross-validation among all nodewise regressions. The subscript "sc" stands for the scaled Lasso. $t(4)/\sqrt{2}$ and Gamma denote the studentized $t(4)$ distribution and the centralized and studentized Gamma(4,1) distribution respectively. S&C denotes the screen and clean procedure. The mean, SD, FP and FN are computed based on $1,000$ simulation runs.

Table 4: Empirical sizes (upper panel) and powers (lower panel) based on the one-step (without sample splitting and screening) and three-step procedures with the non-studentized ("NST") and the studentized ("ST") statistics, where $\mathcal{S}_0 = \{1,2,3\}$ and $\widetilde{\mathcal{S}}_0 = \{1,2,3,\ldots,15\}$, and $p = 500$. Case (i) ((ii)) corresponds to Toeplitz matrix (exchangeable matrix). The nominal levels are 5% and 1%.

|  |  | One-Step Procedure | | | | Three-Step Procedure | | | |
|---|---|---|---|---|---|---|---|---|---|
|  |  | $t(4)/\sqrt{2}$ | | Gamma | | $t(4)/\sqrt{2}$ | | Gamma | |
|  | $\alpha$ | NST | ST | NST | ST | NST | ST | NST | ST |
| (i), $\mathcal{S}_0^c$ | 5% | 0.03 | 0.03 | 0.02 | 0.03 | 0.06 | 0.06 | 0.07 | 0.07 |
|  | 1% | 0.01 | 0.00 | 0.00 | 0.01 | 0.01 | 0.01 | 0.03 | 0.02 |
| (i), $\widetilde{\mathcal{S}}_0^c$ | 5% | 0.02 | 0.01 | 0.02 | 0.02 | 0.05 | 0.04 | 0.06 | 0.05 |
|  | 1% | 0.01 | 0.00 | 0.01 | 0.01 | 0.02 | 0.02 | 0.03 | 0.02 |
| (ii), $\mathcal{S}_0^c$ | 5% | 0.09 | 0.08 | 0.08 | 0.08 | 0.10 | 0.11 | 0.09 | 0.10 |
|  | 1% | 0.02 | 0.03 | 0.03 | 0.03 | 0.04 | 0.04 | 0.03 | 0.03 |
| (i), $\{3\} \cup \mathcal{S}_0^c$ | 5% | 0.66 | 0.67 | 0.61 | 0.62 | 0.74 | 0.72 | 0.73 | 0.72 |
|  | 1% | 0.54 | 0.56 | 0.49 | 0.51 | 0.62 | 0.60 | 0.63 | 0.61 |
| (i), $\{2,3\} \cup \mathcal{S}_0^c$ | 5% | 0.89 | 0.88 | 0.88 | 0.87 | 0.93 | 0.92 | 0.94 | 0.93 |
|  | 1% | 0.80 | 0.80 | 0.79 | 0.78 | 0.87 | 0.85 | 0.87 | 0.85 |
| (i), $\{15\} \cup \widetilde{\mathcal{S}}_0^c$ | 5% | 0.63 | 0.64 | 0.58 | 0.59 | 0.73 | 0.72 | 0.70 | 0.68 |
|  | 1% | 0.53 | 0.54 | 0.48 | 0.50 | 0.63 | 0.61 | 0.60 | 0.57 |
| (i), $\{14,15\} \cup \widetilde{\mathcal{S}}_0^c$ | 5% | 0.86 | 0.86 | 0.82 | 0.82 | 0.93 | 0.92 | 0.93 | 0.92 |
|  | 1% | 0.76 | 0.76 | 0.71 | 0.72 | 0.86 | 0.85 | 0.87 | 0.85 |
| (ii), $\{3\} \cup \mathcal{S}_0^c$ | 5% | 1.00 | 1.00 | 1.00 | 1.00 | 0.99 | 0.99 | 1.00 | 1.00 |
|  | 1% | 1.00 | 1.00 | 1.00 | 1.00 | 0.99 | 0.99 | 1.00 | 1.00 |
| (ii), $\{2,3\} \cup \mathcal{S}_0^c$ | 5% | 1.00 | 1.00 | 1.00 | 1.00 | 1.00 | 1.00 | 1.00 | 1.00 |
|  | 1% | 1.00 | 1.00 | 1.00 | 1.00 | 1.00 | 1.00 | 1.00 | 1.00 |

Note: The tuning parameters $\lambda_j$s in the nodewise Lasso are chosen to be the same via 10-fold cross-validation among all nodewise regressions. $t(4)/\sqrt{2}$ and Gamma denote the studentized $t(4)$ distribution and the centralized and studentized Gamma(4,1) distribution respectively. The sizes and powers are computed based on 1,000 simulation runs.

Table 5: FWER and power of multiple testing based on the step-down method with the non-studentized ("NST") and studentized ("ST") test statistics, and based on the Bonferroni-Holm procedure ("BH"), where $p = 500$, and the nominal level is 5%. Case (i) ((ii)) corresponds to Toeplitz matrix (block diagonal matrix).

|  |  | $s_0 = 3$, (i) | | $s_0 = 15$, (i) | | $s_0 = 3$, (ii) | | $s_0 = 15$, (ii) | |
| --- | --- | --- | --- | --- | --- | --- | --- | --- | --- |
|  |  | FWER | Power | FWER | Power | FWER | Power | FWER | Power |
| $t(4)/\sqrt{2}$ | NST | 0.037 | 0.548 | 0.008 | 0.732 | 0.046 | 0.594 | 0.005 | 0.701 |
|  | ST | 0.028 | 0.534 | 0.009 | 0.722 | 0.048 | 0.560 | 0.005 | 0.685 |
|  | BH | 0.024 | 0.528 | 0.006 | 0.717 | 0.040 | 0.555 | 0.004 | 0.678 |
| Gamma | NST | 0.034 | 0.535 | 0.014 | 0.725 | 0.039 | 0.581 | 0.001 | 0.701 |
|  | ST | 0.033 | 0.513 | 0.016 | 0.714 | 0.046 | 0.549 | 0.004 | 0.685 |
|  | BH | 0.023 | 0.506 | 0.011 | 0.708 | 0.038 | 0.545 | 0.004 | 0.680 |

Note: The tuning parameters $\lambda_j$s in the nodewise Lasso are chosen to be the same via 10-fold cross-validation among all nodewise regressions. $t(4)/\sqrt{2}$ and Gamma denote the studentized $t(4)$ distribution and the centralized and studentized Gamma(4,1) distribution respectively. The FWER and powers are computed based on 1,000 simulation runs.

# Supplement to "Simultaneous Inference for High-dimensional Linear Models"

Xianyang Zhang* and Guang Cheng†

*Texas A&M University and Purdue University*

February 24, 2016

This supplementary material provides proofs of the main results in the paper as well as some additional numerical results.

# 1 Technical details

We first present two lemmas that will be used in the rest proofs. Define $\xi_{ij} = \Theta_j^T \widetilde{X}_i \epsilon_i$. Denote by $c, c', C, C', C_i$ be some generic constants which can be different from line to line.

LEMMA **1.1**. *Under Assumptions 2.1-2.3, we have for any $G \subseteq \{1, 2, \ldots, p\}$,*

$$\sup_{x \in \mathbb{R}} \left| P\left( \max_{j \in G} \sum_{i=1}^n \xi_{ij}/\sqrt{n} \leq x \right) - P\left( \max_{j \in G} \sum_{i=1}^n z_{ij}/\sqrt{n} \leq x \right) \right| \lesssim n^{-c'}, \quad c' > 0,$$

*where $\{z_i = (z_{i1}, \ldots, z_{ip})'\}$ is a sequence of mean zero independent Gaussian vector with $\mathbb{E} z_i z_i' = \Theta_j^T \Sigma \Theta_j \sigma_\epsilon^2$.*

*Assistant Professor, Department of Statistics, Texas A&M University, College Station, TX 77843. E-mail: zhangxiany@stat.tamu.edu.

†Associate Professor, Department of Statistics, Purdue University, West Lafayette, IN 47906. E-mail: chengg@purdue.edu. Tel: +1 (765) 496-9549. Fax: +1 (765) 494-0558. Research Sponsored by NSF CAREER Award DMS-1151692, DMS-1418042, Simons Fellowship in Mathematics, Office of Naval Research (ONR N00014-15-1-2331) and a grant from Indiana Clinical and Translational Sciences Institute. Guang Cheng was on sabbatical at Princeton while part of this work was carried out; he would like to thank the Princeton ORFE department for its hospitality and support.

*Proof of Lemma 1.1.* We apply Corollary 2.1 of Chernozhukov et al. (2013) to the sequence $\{\xi_{ij}\}$ by verifying its Condition (E.1). For the sake of clarity, we state the condition below, i.e.

$$c_0 \leq \mathbb{E}\xi_{ij}^2 \leq C_0, \quad \max_{k=1,2} \mathbb{E}|\xi_{ij}|^{2+k}/B^k + \mathbb{E}\exp(|\xi_{ij}|/B) \leq 4, \tag{1}$$

uniformly over $j$, where $c_0, C_0 > 0$, and $B$ is some large enough constant. In what follows, we consider two cases for $\mathbf{X}$: (i) $\mathbf{X}$ has i.i.d. sub-Gaussian rows; (ii) $\mathbf{X}$ is strongly bounded.

(i) By Assumption 2.2, $\mathbb{E}(\Theta_j^T \widetilde{X}_i)^2 = \Theta_j^T \Sigma \Theta_j = \theta_{jj} := 1/\tau_j^2$, and $1/c < \Lambda_{\min}^2 \leq \tau_j^2 \leq \Sigma_{j,j} = C$, for some constants $c, C > 0$. Recall that $\Lambda_{\min}^2$ is the minimal eigenvalue of $\Sigma$. Thus we have $c_1 \sigma_\epsilon^2 \leq \mathbb{E}\xi_{ij}^2 \leq C_1 \sigma_\epsilon^2$. By the independence between $\{\widetilde{X}_i\}$ and $\{\epsilon_i\}$, we have for large enough $C$ and uniformly for all $j$,

$$\mathbb{E}\exp(|\xi_{ij}|/C) = 1 + \sum_{k=1}^{+\infty} \frac{\mathbb{E}|\xi_{ij}|^k}{C^k k!} = 1 + \sum_{k=1}^{+\infty} \frac{\mathbb{E}|\Theta_j^T \widetilde{X}_i|^k \mathbb{E}|\epsilon_i|^k}{C^k k!}$$

$$\leq 1 + \sum_{k=1}^{+\infty} \frac{k^k}{(C')^k k!} \leq 1 + \sum_{k=1}^{+\infty} (e/C')^k < \infty,$$

where we have used the fact that $k! \geq (k/e)^k$, $\|\Theta_j\|_2 \lesssim \Lambda_{\min}^{-1} = O(1)$ (because $\|\Theta_j\|_2^2 \Lambda_{\min}^2 \leq c$) and $\mathbb{E}|X|^k \leq (C'')^k k^{k/2}$ with $C''$ being some positive constant for sub-Gaussian variable $X$. Thus we have $\max_{k=1,2} \mathbb{E}|\xi_{ij}|^{2+k}/B^k + \mathbb{E}\exp(|\xi_{ij}|/B) \leq 4$ uniformly for some large enough constant $B$.

(ii) In the strongly bounded case, using the fact that $\|\Theta_j\|_2^2 \lesssim \Lambda_{\min}^{-2} = O(1)$ and $\|\Theta_j\|_1 \leq \sqrt{s_j}\|\Theta_j\|_2$, we have $|\Theta_j^T \widetilde{X}_i| \leq \|\Theta_j\|_1 \|\widetilde{X}_i\|_\infty \leq K_n \sqrt{s_j} \|\Theta_j\|_2$. It is straightforward to verify that $\max_{k=1,2} \mathbb{E}|\xi_{ij}|^{2+k}/B_n^k + \mathbb{E}\exp(|\xi_{ij}|/B_n) \leq 4$ uniformly with some $B_n \asymp K_n \max_j \sqrt{s_j}$ and $B_n^2 (\log(pn))^7/n \leq C_2 n^{-c_2}$ under part (ii) of Assumption 2.3. $\diamondsuit$

REMARK 1.1. The conclusion in Lemma 1.1 still holds if we assume that (i) $\max_{i,j} |X_{ij}| \leq K_n$ with $\max_{1 \leq j \leq p} s_j^2 K_n^4 (\log(pn))^7/n \leq C_1 n^{-c_1}$ for some constants $c_1, C_1 > 0$; and (ii) $\{\epsilon_i\}$ are i.i.d with with $\mathbb{E}|\epsilon_i|^4 < \infty$ and $c' < \sigma_\epsilon^2$ for $c' > 0$.

Next we quantify the effect by replacing $\xi_i$ with $\widehat{\xi}_i$.

LEMMA **1.2**. *Suppose Assumptions 2.1-2.3 hold. Assume $\max_j K_0^2 s_j^2 (\log(pn))^3 (\log(n))^2/n = o(1)$. Recall that $K_0 = 1$ in the sub-Gaussian case and $K_0 = K_n$ in the strongly bounded case. Then with $\lambda_j \asymp K_0 \sqrt{\log(p)/n}$ uniformly for $j$, there exist $\zeta_1, \zeta_2 > 0$ such that*

$$P\left(\max_{1 \le j \le p} \left| \sum_{i=1}^n \widehat{\xi}_{ij}/\sqrt{n} - \sum_{i=1}^n \xi_{ij}/\sqrt{n} \right| \ge \zeta_1 \right) < \zeta_2,$$

*where $\zeta_1 \sqrt{1 \vee \log(p/\zeta_1)} = o(1)$ and $\zeta_2 = o(1)$.*

*Proof of Lemma 1.2.* Let $\widetilde{K}_0 = \log(np)\log(n)$ in the sub-Gaussian case and $\widetilde{K}_0 = K_n \log(n)$ in the strongly bounded case. Using Lemma A.1 in Chernozhukov et al. (2013), we deduce that

$$\begin{aligned}
\mathbb{E}\left\{ \max_{1 \le j \le p} \left| \sum_{i=1}^n X_{ij}\epsilon_i/n \right| \right\} &\lesssim \sigma_\epsilon \sqrt{\max_j \Sigma_{j,j}} \sqrt{\log(p)/n} + \sqrt{\mathbb{E} \max_{i,j} |X_{ij}\epsilon_i|^2 \log(p)/n} \\
&\lesssim \sqrt{\log(p)/n} + \sqrt{\mathbb{E} \max_{i,j} X_{ij}^2} \sqrt{\mathbb{E} \max_i \epsilon_i^2 \log(p)/n} \\
&\lesssim \sqrt{\log(p)/n} + \widetilde{K}_0 \log(p)/n,
\end{aligned}$$

where we have used the fact that $\sqrt{\mathbb{E} \max_i \epsilon_i^2} \lesssim \log(n) \max_{1 \le i \le n} ||\epsilon_i||_{\psi_1} \lesssim \log n$ with $\psi_1(x) = \exp(x) - 1$ and $||\cdot||_{\psi_1}$ being the corresponding Orlicz norm, and similar result for $\sqrt{\mathbb{E} \max_{i,j} X_{i,j}^2}$ (see Lemma 2.2.2 in van der Vaart and Wellner 1996). Because $||\widehat{\Theta}_j - \Theta_j||_1 = O_P(K_0 s_j \sqrt{\log(p)/n})$ uniformly for $j$, we obtain,

$$\begin{aligned}
\left| \sum_{i=1}^n \widehat{\xi}_{ij}/\sqrt{n} - \sum_{i=1}^n \xi_{ij}/\sqrt{n} \right| &= \left| (\widehat{\Theta}_j^T - \Theta_j^T) \sum_{i=1}^n \widetilde{X}_i \epsilon_i / \sqrt{n} \right| \le ||\widehat{\Theta}_j - \Theta_j||_1 \left\| \sum_{i=1}^n \widetilde{X}_i \epsilon_i / \sqrt{n} \right\|_\infty \\
&= O_P\left( K_0 s_j \sqrt{\log(p)/n} \left\| \sum_{i=1}^n \widetilde{X}_i \epsilon_i / \sqrt{n} \right\|_\infty \right) \\
&= O_P\left( K_0 s_j \log(p)/\sqrt{n} + \sqrt{n} K_0 \widetilde{K}_0 s_j (\log(p)/n)^{3/2} \right) \\
&\le O_P\left( \max_j s_j K_0 \log(p)/\sqrt{n} \right),
\end{aligned}$$

uniformly for all $j$. Choosing $\zeta_1$ such that $\max_j K_0 s_j \log(p)/(\sqrt{n}\zeta_1) = o(1)$ and $\zeta_1 \sqrt{1 \vee \log(p/\zeta_1)} =$

$o(1)$ (e.g. $\zeta_1^2 = O(\max_j K_0 s_j \sqrt{\log(p)/n})$), we deduce that

$$P\left(\max_{1\leq j\leq p}\left|\sum_{i=1}^n \widehat{\xi}_{ij}/\sqrt{n} - \sum_{i=1}^n \xi_{ij}/\sqrt{n}\right| \geq \zeta_1\right) < \zeta_2, \quad \zeta_2 = o(1).$$

◇

REMARK 1.2. With a more delicate analysis, one can specify the order of $\zeta_2$ in Lemma 1.2; see e.g., Theorem 6.1 and Lemma 6.2 of Bühlmann and van de Geer (2011).

*Proof of Theorem 2.2.* Without loss of generality, we set $G = \{1, 2, \ldots, p\}$. Define

$$T_G = \max_{j\in G} \sqrt{n}(\breve{\beta}_j - \beta_j^0), \quad T_{0,G} = \max_{j\in G} \sum_{i=1}^n \xi_{ij}/\sqrt{n}.$$

Let $\pi(v) = C_2 v^{1/3}(1 \vee \log(p/v))^{2/3}$ with $C_2 > 0$, and

$$\Gamma = \max_{1\leq j,k\leq p} |\widehat{\sigma}_\epsilon^2 \widehat{\Theta}_j^T \widehat{\Sigma} \widehat{\Theta}_k - \sigma_\epsilon^2 \Theta_j^T \Sigma \Theta_k|, \quad \widehat{\Sigma} = \mathbf{X}^T\mathbf{X}/n.$$

Notice that

$$|T_G - T_{0,G}| \leq \max_{1\leq j\leq p}\left|\sum_{i=1}^n \widehat{\xi}_{ij}/\sqrt{n} - \sum_{i=1}^n \xi_{ij}/\sqrt{n}\right| + ||\Delta||_\infty.$$

By similar arguments in the proof of Theorem 2.4 of van de Geer et al. (2014) and the results in Theorem 2.1, we have

$$||\Delta||_\infty \leq ||\widehat{\beta} - \beta^0||_1 \max_j \sqrt{n}\lambda_j/\widehat{\tau}_j^2 = O_P(K_0\sqrt{\log(p)}||\widehat{\beta} - \beta^0||_1) = O_P(K_0^2 s_0 \log(p)/\sqrt{n}),$$

where we use the fact that $\max_j \lambda_j/\widehat{\tau}_j^2 = O_P(K_0\sqrt{\log(p)/n})$ and $||\widehat{\beta} - \beta^0||_1 = O_P(s_0\lambda)$ with $\lambda = O(K_0\sqrt{\log(p)/n})$. Thus by Lemma 1.2 and the assumption that $K_0^4 s_0^2(\log(p))^3/n = o(1)$, we have

$$P(|T_G - T_{0,G}| > \zeta_1) < \zeta_2,$$

for $\zeta_1\sqrt{1 \vee \log(p/\zeta_1)} = o(1)$ and $\zeta_2 = o(1)$.

Let $c_{z,G}(\alpha) = \inf\{t \in \mathbb{R} : P(\max_{j\in G}\sum_{i=1}^n z_{ij}/\sqrt{n} \leq t) \geq 1 - \alpha\}$, where the sequence $\{z_{ij}\}$ is

defined in Lemma 1.1. Following the arguments in the proof of Lemma 3.2 in Chernozhukov et al. (2013), we have

$$P(c_G(\alpha) \leq c_{z,G}(\alpha + \pi(v))) \geq 1 - P(\Gamma > v), \tag{2}$$

$$P(c_{z,G}(\alpha) \leq c_G(\alpha + \pi(v))) \geq 1 - P(\Gamma > v). \tag{3}$$

By Lemma 1.1, (2) and (3), we have for every $v > 0$,

$$\sup_{\alpha \in (0,1)} |P(T_{0,G} > c_G(\alpha)) - \alpha| \lesssim \sup_{\alpha \in (0,1)} \left| P\left( \max_{j \in G} \sum_{i=1}^n z_{ij}/\sqrt{n} > c_G(\alpha) \right) - \alpha \right| + n^{-c'}$$

$$\lesssim \pi(v) + P(\Gamma > v) + n^{-c'}.$$

Moreover, by the arguments in the proof of Theorem 3.2 in Chernozhukov et al. (2013), we have

$$\sup_{\alpha \in (0,1)} |P(T_G > c_G(\alpha)) - \alpha| \lesssim \pi(v) + P(\Gamma > v) + n^{-c'} + \zeta_1\sqrt{1 \vee \log(p/\zeta_1)} + \zeta_2.$$

By Lemma 5.3 and Lemma 5.4 of van de Geer et al. (2014), we have

$$\max_{1 \leq j,k \leq p} |\widehat{\Theta}_j^T \widehat{\Sigma} \widehat{\Theta}_k - \Theta_j^T \Sigma \Theta_k| = O_P(\max_j \lambda_j \sqrt{s_j}).$$

Since $|\Theta_j^T \Sigma \Theta_k| \leq 1/(\tau_j \tau_k) = O(1)$ uniformly for $1 \leq j, k \leq p$, we have

$$\Gamma = O_P\left( |\widehat{\sigma}_\epsilon^2 - \sigma_\epsilon^2| + \max_j \lambda_j \sqrt{s_j} \right).$$

Under Assumption 2.4, choosing $v = 1/(\alpha_n (\log(p))^2)$, we deduce that

$$\sup_{\alpha \in (0,1)} \left| P(\max_{1 \leq j \leq p} \sqrt{n}(\breve{\beta}_j - \beta_j^0) > c_G(\alpha)) - \alpha \right| = o(1),$$

which completes the proof. $\diamondsuit$

*Proof of Theorem 2.3.* From the arguments in the proof of Theorem 2.2, we have

$$\Gamma = \max_{1\leq j,k \leq p} |\widehat{\sigma}_\epsilon^2 \widehat{\Theta}_j^T \widehat{\Sigma} \widehat{\Theta}_k - \sigma_\epsilon^2 \Theta_j^T \Sigma \Theta_k| = O_P\left(|\widehat{\sigma}_\epsilon^2 - \sigma_\epsilon^2| + \max_j \lambda_j \sqrt{s_j}\right),$$

which implies that $\max_{1\leq j \leq p} |\widehat{\omega}_{jj} - \omega_{jj}| = O_P\left(|\widehat{\sigma}_\epsilon^2 - \sigma_\epsilon^2| + \max_j \lambda_j \sqrt{s_j}\right)$ with $\omega_{jj} = \sigma_\epsilon^2 \theta_{jj}$. We then have

$$P(\omega_{jj}/2 < \widehat{\omega}_{jj} < 2\omega_{jj} \text{ for all } 1 \leq j \leq p) \to 1. \tag{4}$$

The fact that $1/c < \Lambda_{\min}^2 \leq \tau_j^2 = 1/\theta_{jj} \leq \Sigma_{j,j} = C$ implies that $\omega_{jj}$ is uniformly bounded away from zero and infinity.

Define $\bar{T}_G = \max_{j \in G} \sqrt{n}(\breve{\beta}_j - \beta_j^0)/\sqrt{\widehat{\omega}_{jj}}$ and $\bar{T}_{0,G} = \max_{j \in G} \sum_{i=1}^n \xi_{ij}/\sqrt{n\omega_{jj}}$. Denote by $\Delta = (\Delta_1, \ldots, \Delta_p)^T$ and $\bar{\Delta} = (\bar{\Delta}_1, \ldots, \bar{\Delta}_p)^T$ with $\bar{\Delta}_j = \Delta_j/\sqrt{\widehat{\omega}_{jj}}$. Then we have

$$|\bar{T}_G - \bar{T}_{0,G}|$$

$$\leq \max_{1\leq j \leq p} \left| \sum_{i=1}^n \widehat{\xi}_{ij}/\sqrt{n\widehat{\omega}_{jj}} - \sum_{i=1}^n \xi_{ij}/\sqrt{n\omega_{jj}} \right| + ||\bar{\Delta}||_\infty$$

$$\leq \max_{1\leq j \leq p} \left| \sum_{i=1}^n \widehat{\xi}_{ij}/\sqrt{n\widehat{\omega}_{jj}} - \sum_{i=1}^n \widehat{\xi}_{ij}/\sqrt{n\omega_{jj}} \right| + \max_{1\leq j \leq p} \left| \sum_{i=1}^n \widehat{\xi}_{ij}/\sqrt{n\omega_{jj}} - \sum_{i=1}^n \xi_{ij}/\sqrt{n\omega_{jj}} \right| + ||\bar{\Delta}||_\infty$$

$$\leq C' \max_{1\leq j \leq p} \left| \sum_{i=1}^n \widehat{\xi}_{ij}/\sqrt{n} \right| \max_{1\leq j \leq p} \left| \sqrt{\omega_{jj}/\widehat{\omega}_{jj}} - 1 \right| + C'' \max_{1\leq j \leq p} \left| \sum_{i=1}^n (\widehat{\xi}_{ij} - \xi_{ij})/\sqrt{n} \right| + ||\bar{\Delta}||_\infty,$$

$$= I_1 + I_2 + I_3,$$

where $C', C'' > 0$.

On the event $\omega_{jj}/2 < \widehat{\omega}_{jj} < 2\omega_{jj}$ for all $1 \leq j \leq p$,

$$\max_{1 \leq j \leq p} \left| \sqrt{\omega_{jj}/\widehat{\omega}_{jj}} - 1 \right| \leq \max_{1 \leq j \leq p} |\sqrt{\omega_{jj}} - \sqrt{\widehat{\omega}_{jj}}| \max_{1 \leq j \leq p} \sqrt{2/\omega_{jj}}$$

$$\leq \max_{1 \leq j \leq p} \left| \frac{\omega_{jj} - \widehat{\omega}_{jj}}{\sqrt{\omega_{jj}} + \sqrt{\widehat{\omega}_{jj}}} \right| \max_{1 \leq j \leq p} \sqrt{2/\omega_{jj}}$$

$$\leq \max_{1 \leq j \leq p} |\omega_{jj} - \widehat{\omega}_{jj}| \max_{1 \leq j \leq p} 1/\omega_{jj}$$

$$= O_P \left( |\widehat{\sigma}_\epsilon^2 - \sigma_\epsilon^2| + \max_j \lambda_j \sqrt{s_j} \right).$$

On the other hand,

$$\max_{1 \leq j \leq p} \left| \sum_{i=1}^n \widehat{\xi}_{ij}/\sqrt{n} \right| \leq \max_{1 \leq j \leq p} \left| \sum_{i=1}^n (\widehat{\xi}_{ij} - \xi_{ij})/\sqrt{n} \right| + \max_{1 \leq j \leq p} \left| \sum_{i=1}^n \xi_{ij}/\sqrt{n} \right|$$

$$= O_P(\sqrt{\log(p)} + \max_j \sqrt{s_j} \widetilde{K}_0 \log(p)/\sqrt{n}) = O_P(\sqrt{\log p}),$$

where $\widetilde{K}_0 = \log(np)\log(n)$ in the sub-Gaussian case and $\widetilde{K}_0 = K_n \log(n)$ in the strongly bounded case. Therefore, on the above event, $I_1 \leq O_P \left( \sqrt{\log(p)} |\widehat{\sigma}_\epsilon^2 - \sigma_\epsilon^2| + \sqrt{\log(p)} \max_j \lambda_j \sqrt{s_j} \right)$. Under Assumption 2.4, we can find $\zeta_1'$ such that $P(I_1 > \zeta_1') = o(1)$ and $\zeta_1' \sqrt{1 \vee \log(p/\zeta_1')} = o(1)$. Using the fact that $||\Delta||_\infty \leq O_P(K_0^2 s_0 \log(p)/\sqrt{n})$, we can prove the same result for $||\bar{\Delta}||_\infty$ conditional on the event $\{\omega_{jj}/2 < \widehat{\omega}_{jj} < 2\omega_{jj}$ for all $1 \leq j \leq p\}$. Thus by Lemma 1.2 and (4), we have

$$P(|\bar{T}_G - \bar{T}_{0,G}| > \zeta_1) \leq P(I_1 + I_2 + I_3 > \zeta_1) < \zeta_2,$$

for $\zeta_1 \sqrt{1 \vee \log(p/\zeta_1)} = o(1)$ and $\zeta_2 = o(1)$.

Let $\bar{\Gamma} = \max_{1 \leq j,k \leq p} |\widehat{\sigma}_\epsilon^2 \widehat{\Theta}_j^T \widehat{\Sigma} \widehat{\Theta}_k / \sqrt{\widehat{\omega}_{jj} \widehat{\omega}_{kk}} - \sigma_\epsilon^2 \Theta_j^T \Sigma \Theta_k / \sqrt{\omega_{jj} \omega_{kk}}|$. Note that

$$|\sqrt{\omega_{jj} \omega_{kk}} - \sqrt{\widehat{\omega}_{jj} \widehat{\omega}_{kk}}| = \frac{|\omega_{jj} \omega_{kk} - \widehat{\omega}_{jj} \widehat{\omega}_{kk}|}{\sqrt{\omega_{jj} \omega_{kk}} + \sqrt{\widehat{\omega}_{jj} \widehat{\omega}_{kk}}}.$$

On the event $\omega_{jj}/2 < \widehat{\omega}_{jj} < 2\omega_{jj}$ for all $1 \leq j \leq p$, we have

$$\frac{|\omega_{jj}\omega_{kk} - \widehat{\omega}_{jj}\widehat{\omega}_{kk}|}{\sqrt{\omega_{jj}\omega_{kk}} + \sqrt{\widehat{\omega}_{jj}\widehat{\omega}_{kk}}} \leq \frac{|\omega_{jj}\omega_{kk} - \widehat{\omega}_{jj}\widehat{\omega}_{kk}|}{\sqrt{\omega_{jj}\omega_{kk}} + \sqrt{\omega_{jj}\omega_{kk}/4}} \leq (2/3)|\omega_{jj}\omega_{kk} - \widehat{\omega}_{jj}\widehat{\omega}_{kk}| \max_{1 \leq j \leq p} 1/\omega_{jj},$$

which implies that

$$\max_{1 \leq j,k \leq p} |\sqrt{\omega_{jj}\omega_{kk}/\widehat{\omega}_{jj}\widehat{\omega}_{kk}} - 1| \leq \max_{1 \leq j,k \leq p} |\sqrt{\omega_{jj}\omega_{kk}} - \sqrt{\widehat{\omega}_{jj}\widehat{\omega}_{kk}}| \max_{1 \leq j \leq p} 2/\omega_{jj}$$

$$\leq (4/3) \max_{1 \leq j,k \leq p} |\omega_{jj}\omega_{kk} - \widehat{\omega}_{jj}\widehat{\omega}_{kk}| \max_{1 \leq j \leq p} 1/\omega_{jj}^2$$

$$= O_P\left(|\widehat{\sigma}_\epsilon^2 - \sigma_\epsilon^2| + \max_j \lambda_j \sqrt{s_j}\right).$$

Using similar arguments above, we can show that $P(\bar{\Gamma} > v) = o(1)$ for $v = 1/(\alpha_n(\log(p))^2)$. The rest of the proofs is similar to those in the proof of Theorem 2.2. We skip the details $\diamondsuit$

*Proof of Theorem 2.4.* Define $\widetilde{T}_G = \max_{j \in G} |\sqrt{n}(\breve{\beta}_j - \beta_j^0)/\sqrt{\widehat{\omega}_{jj}}|$ and $\widetilde{T}_{0,G} = \max_{j \in G} \sum_{i=1}^n |\xi_{ij}/\sqrt{n\omega_{jj}}|$. Under the assumptions in Theorem 2.3, we can show that $P(|\widetilde{T}_G - \widetilde{T}_{0,G}| > \zeta_1) < \zeta_2$ for $\zeta_1\sqrt{1 \vee \log(p/\zeta_1)} = o(1)$ and $\zeta_2 = o(1)$. In another word, the distribution of $\max_{j \in G} \sqrt{n}|\breve{\beta}_j - \beta_0|/\sqrt{\widehat{\omega}_{jj}}$ can be approximated by $\max_{j \in G} |Z_j|$ with $Z = (Z_1, \ldots, Z_p) \sim^d N(0, \widetilde{\Theta})$. Under Assumption 2.5, by Lemma 6 of Cai et al. (2014), we have for any $x \in \mathbb{R}$ and as $|G| \to +\infty$,

$$P\left(\max_{j \in G} |Z_i|^2 - 2\log(|G|) + \log\log(|G|) \leq x\right) \to F(x) := \exp\left\{-\frac{1}{\sqrt{\pi}}\exp\left(-\frac{x}{2}\right)\right\}.$$

It implies that

$$P\left(\max_{j \in G} n|\breve{\beta}_j - \beta_j^0|^2/\widehat{\omega}_{jj} \leq 2\log(|G|) - \log\log(|G|)/2\right) \to 1. \quad (5)$$

The bootstrap consistency result implies that

$$|(\bar{c}_G^*(\alpha))^2 - 2\log(|G|) + \log\log(|G|) - q_\alpha| = o_P(1), \quad (6)$$

where $q_\alpha$ is the $100(1-\alpha)$th quantile of $F(x)$. Consider any $j^* \in G$ such that $|\breve{\beta}_{j^*} - \beta_{j^*}^0|/\sqrt{\omega_{j^*j^*}} >$

$(\sqrt{2} + \varepsilon_0)\sqrt{(\log |G|)/n}$. Using the inequality $2a_1 a_2 \leq \delta^{-1} a_1^2 + \delta a_2^2$ for any $\delta > 0$, we have

$$n|\widetilde{\beta}_{j^*} - \beta_{j^*}^0|^2/\widehat{\omega}_{j^* j^*} \leq (1+\delta^{-1}) n|\breve{\beta}_{j^*} - \beta_{j^*}^0|^2/\widehat{\omega}_{j^* j^*} + (1+\delta) n|\breve{\beta}_{j^*} - \widetilde{\beta}_{j^*}|^2/\widehat{\omega}_{j^* j^*}, \qquad (7)$$

where $n|\breve{\beta}_{j^*} - \beta_{j^*}^0|^2/\widehat{\omega}_{j^* j^*} = o_p(\log |G|)$ as $j^*$ is fixed and $|G|$ grows. From the proof of Theorem 2.3, we know the difference between $n|\widetilde{\beta}_{j^*} - \beta_{j^*}^0|^2/\widehat{\omega}_{j^* j^*}$ and $n|\widetilde{\beta}_{j^*} - \beta_{j^*}^0|^2/\omega_{j^* j^*}$ is asymptotically negligible. Thus by (7) and the fact that $\beta^0 \in \mathcal{U}_G(\sqrt{2} + \varepsilon_0)$, we have,

$$\max_{j \in G} n|\breve{\beta}_j - \widetilde{\beta}_j|^2/\widehat{\omega}_{jj} \geq \frac{1}{1+\delta} \left\{ (\sqrt{2} + \varepsilon_0)^2 (\log |G|) - o_p(\log |G|) \right\}. \qquad (8)$$

The conclusion thus follows from (8) and (6) provided that $\delta$ is small enough. $\diamondsuit$

*Proof of Proposition 3.1.* Similar to the proof of Theorem 2.4, the distribution of $\max_{1 \leq j \leq p} \sqrt{n}|\breve{\beta}_j - \beta_0|/\sqrt{\widehat{\omega}_{jj}}$ can be approximated by $\max_{1 \leq j \leq p} |Z_j|$ with $Z = (Z_1, \ldots, Z_p) \overset{d}{\sim} N(0, \widetilde{\Theta})$. Under Assumption 2.5, by Lemma 6 of Cai et al. (2014), we have for any $x \in \mathbb{R}$ and as $p \to +\infty$,

$$P\left( \max_{1 \leq i \leq p} |Z_i|^2 - 2\log(p) + \log\log(p) \leq x \right) \to \exp\left\{ -\frac{1}{\sqrt{\pi}} \exp\left( -\frac{x}{2} \right) \right\}.$$

It implies that

$$P\left( \max_{j \in \mathcal{S}_0^c} n|\breve{\beta}_j|^2/\widehat{\omega}_{jj} \leq 2\log(p) - \log\log(p)/2 \right) \to 1. \qquad (9)$$

On the other hand, we note that

$$\min_{j \in \mathcal{S}_0} n|\beta_j^0|^2/\widehat{\omega}_{jj} \leq 2 \max_{j \in \mathcal{S}_0} n|\breve{\beta}_j - \beta_j^0|^2/\widehat{\omega}_{jj} + 2 \min_{j \in \mathcal{S}_0} n|\breve{\beta}_j|^2/\widehat{\omega}_{jj}$$

Because the difference between $\min_{j \in \mathcal{S}_0} n|\beta_j^0|^2/\widehat{\omega}_{jj}$ and $\min_{j \in \mathcal{S}_0} n|\beta_j^0|^2/\omega_{jj}$ is asymptotically negligible, and $P(2 \max_{j \in \mathcal{S}_0} n|\breve{\beta}_j - \beta_j^0|^2/\widehat{\omega}_{jj} \leq 4\log(p) - \log\log(p)) \to 1$, we obtain

$$\begin{aligned} &P\left( \min_{j \in \mathcal{S}_0} n|\breve{\beta}_j|^2/\widehat{\omega}_{jj} > 2\log p \right) \\ &\geq P\left( 2 \min_{j \in \mathcal{S}_0} n|\breve{\beta}_j|^2/\widehat{\omega}_{jj} + 4\log(p) - \log\log(p) > 8\log(p) \right) \to 1. \end{aligned} \qquad (10)$$

Hence, (16) follows from (9) and (10).

We next prove the optimality of $\tau^* = 2$, i.e., (17). For large enough $p$, we can choose a set $G^*$ such that $\beta_j = 0$ for $j \in G^*$, and $|G^*| = \lfloor p^{\tau_2} \rfloor$ with $\tau/2 < \tau_2 < 1$. Following the above arguments, we know that the distribution of $\max_{j \in G^*} \sqrt{n}|\breve{\beta}_j - \beta_j^0|/\sqrt{\widehat{\omega}_{jj}}$ can be approximated by $\max_{j \in G^*} |Z_j|$ with $Z = (Z_1, \ldots, Z_p) \sim^d N(0, \widetilde{\Theta})$. Then we have

$$P\left(\max_{j \in G^*} n|\breve{\beta}_j|^2/\widehat{\omega}_{jj} \geq c \log(p)\right) \to 1,$$

where $\tau < c < 2\tau_2 < 2$. The conclusion thus follows immediately. $\diamondsuit$

*Proof of Theorem 4.1.* For simplicity, we only prove the result for the one-sided case (the arguments below can be easily modified for the two-sided case). Define $T_G = \max_{j \in G} \sqrt{n}(\breve{\beta}_j - \beta_j^0)$ and $T_{0,G} = \max_{j \in G} \sum_{i=1}^n \xi_{ij}/\sqrt{n}$. Let $\widetilde{c}_G(\alpha)$ be the bootstrap critical value for the one-sided test at level $\alpha$. We first show that there exist $\zeta_1, \zeta_2 > 0$ such that

$$P\left(|T_G - T_{0,G}| \geq \zeta_1\right) < \zeta_2, \tag{11}$$

where $\zeta_1 \sqrt{1 \vee \log(p/\zeta_1)} = o(1)$ and $\zeta_2 = o(1)$. Notice that

$$|T_G - T_{0,G}| \leq \max_{j \in G} \sqrt{n}|(\Theta_j^T - \widehat{\Theta}_j^T)\mathbb{E}_n \dot{L}_{\beta_0}| + ||\Delta||_\infty + \sqrt{n}||\widehat{\Theta}\mathcal{R}||.$$

Under the Lipschitz continuity in Assumption 4.1, we have

$$\widehat{\Theta}_j^T \mathbb{E}_n \dot{L}_{\widehat{\beta}} = \widehat{\Theta}_j^T \mathbb{E}_n \dot{L}_{\beta_0} + \widehat{\Theta}_j^T \mathbb{E}_n \ddot{L}_{\widehat{\beta}}(\widehat{\beta} - \beta^0) + \mathcal{R}_j,$$

where $\mathcal{R}_j = \widehat{\Theta}_j^T \mathcal{R} \leq \max_i |\widehat{\Theta}_j^T x_i| \cdot ||\mathbf{X}(\widehat{\beta} - \beta^0)||_2^2/n = O_P(K_n s_0 \lambda^2)$ (see the proof of Theorem 3.1 in van de Geer et al. 2014). It thus implies that $\sqrt{n}||\widehat{\Theta}\mathcal{R}||_\infty = O_P(\sqrt{n} K_n s_0 \lambda^2)$. By Assumptions 4.3-4.4, we have

$$||\Delta||_\infty = ||\sqrt{n}(\widehat{\Theta}\widehat{\Sigma} - I)(\widehat{\beta} - \beta_0)||_\infty \leq ||\widehat{\Theta}\widehat{\Sigma} - I||_\infty \sqrt{n}||\widehat{\beta} - \beta_0||_1 = O_P(\sqrt{n}\lambda\lambda_* s_0)$$

Following the arguments in the proof of Lemma 1.2, it can be shown that under Assumption 4.5

$$\max_{j\in G}\sqrt{n}|(\Theta_j^T - \widehat{\Theta}_j^T)\mathbb{E}_n \dot{L}_{\beta_0}| = \max_{j\in G}\left|\sum_{i=1}^n \widehat{\xi}_{ij}/\sqrt{n} - \sum_{i=1}^n \xi_{ij}/\sqrt{n}\right|$$
$$= O_P(K_n \max_j s_j \log(p)/\sqrt{n}) + O_P\left(K_n^2 s_0 \left(\lambda^2 \sqrt{n} \vee \lambda\sqrt{\log(p)}\right)\right)$$

Thus (11) follows from a proper choice of $\zeta_1$.

By Lemma 1.1, we have

$$\sup_{x\in\mathbb{R}}\left|P\left(\max_{j\in G}\sum_{i=1}^n \xi_{ij}/\sqrt{n} \leq x\right) - P\left(\max_{j\in G}\sum_{i=1}^n z_{ij}/\sqrt{n} \leq x\right)\right| \lesssim n^{-c'}, \quad c' > 0,$$

where $\{z_i = (z_{i1}, \ldots, z_{ip})'\}$ is a sequence of mean zero independent Gaussian vector with $\mathbb{E}z_i z_i' = \Theta_j^T \Sigma_{\beta_0} \Theta_j$. By the arguments in the proof of Theorem 3.2 in Chernozhukov et al. (2013), we have

$$\sup_{\alpha\in(0,1)} |P(T_G > \widetilde{c}_G^*(\alpha)) - \alpha| \lesssim \pi(v) + P(\widetilde{\Gamma} > v) + n^{-c'} + \zeta_1\sqrt{1 \vee \log(p/\zeta_1)} + \zeta_2, \quad (12)$$

where $\pi(v) = C_2 v^{1/3}(1 \vee \log(p/v))^{2/3}$. The conclusion follows by choosing $v = 1/(\alpha_n(\log(p))^2)$ in (12). $\diamondsuit$

## 2 Additional numerical results

We consider the linear models where the rows of **X** are fixed i.i.d realizations from $N_p(0, \Sigma)$ with $\Sigma = (\Sigma_{i,j})_{i,j=1}^p$ under two scenarios: (i) Toeplitz: $\Sigma_{i,j} = 0.9^{|i-j|}$; (ii) Exchangeable/Compound symmetric: $\Sigma_{i,i} = 1$ and $\Sigma_{i,j} = 0.8$ for $i \neq j$. The active set is $\mathcal{S}_0 = \{1, 2, \ldots, s_0\}$ with $s_0 = 3$ or 15. To obtain the main Lasso estimator, we implemented the scaled Lasso with the tuning parameter $\lambda_0 = \sqrt{2}\tilde{L}_n(k_0/p)$ with $\tilde{L}_n(t) = n^{-1/2}\Phi^{-1}(1-t)$, where $\Phi$ is the cumulative distribution function for $N(0, 1)$, and $k_0$ is the solution to $k = \tilde{L}_1^4(k/p) + 2\tilde{L}_1^2(k/p)$. We estimate the noise level $\sigma^2$ using the modified variance estimator.

## 2.1 Modified variance estimator

Figure S.1 provides boxplots of $\widehat{\sigma}/\sigma$ for the variance estimator delivered by the scaled Lasso (denoted by "SLasso") and for the modified variance estimator in (24) of the paper (denoted by"SLasso*"). Clearly, the modified variance estimator corrects the noise underestimation issue and thus is preferable.

## 2.2 Impact of the remainder term

We discuss the impact of the (normalized) remainder term $\Delta$ on the coverage accuracy. Recall the linear expansion $\sqrt{n}(\widecheck{\beta} - \beta^0) = \widehat{\Theta}\mathbf{X}^T\epsilon/\sqrt{n} + \Delta$, where $\Delta = (\Delta_1, \ldots, \Delta_p)^T = -\sqrt{n}(\widehat{\Theta}\widehat{\Sigma} - I)(\widehat{\beta} - \beta^0)$ with $\widehat{\Sigma}$ being the Gram matrix and $\widehat{\beta}$ being the Lasso estimator. The studentized maximum type test statistic can be written as

$$\max_{1 \leq j \leq p} \frac{\sqrt{n}|\widecheck{\beta}_j - \beta_j^0|}{\sqrt{\widehat{\omega}_{jj}}} = \max_{1 \leq j \leq p} \left| \frac{\sum_{i=1}^n \widehat{\xi}_{ij}}{\sqrt{n\widehat{\omega}_{jj}}} + \frac{\Delta_j}{\sqrt{\widehat{\omega}_{jj}}} \right|. \tag{13}$$

Thus the coverage accuracy can be greatly affected by the term $\Delta_j^* := \frac{\Delta_j}{\sqrt{\widehat{\omega}_{jj}}}$. Note that this (normalized) remainder term is determined by $\widehat{\Theta}$. We now consider three different methods in estimating $\Theta$ : (i) nodewise Lasso with $\lambda_j$s chosen by 10-fold cross validation; (ii) nodewise Lasso with $\lambda_j = 0.01$; (iii) the method in Javanmard and Montanari (2014) with the tuning parameters chosen automatically by their algorithm. To empirically evaluate $\Delta^* := (\Delta_1^*, \ldots, \Delta_p^*)^T$, we consider the linear models with $t(4)/\sqrt{2}$ errors, $n = 100$ and $p = 500$. Define $\Delta_{ac}^* = (\Delta_j^*)_{j \in \mathcal{S}_0}$ and $\Delta_{in}^* = (\Delta_j^*)_{j \in \mathcal{S}_0^c}$. Figure S.2 presents the boxplots for $||\Delta_{ac}^*||_\infty$ and $||\Delta_{in}^*||_\infty$. The nodewise Lasso clearly outperforms the method in Javanmard and Montanari (2014), and the choice of $\lambda_j = 0.01$ yields the smallest $||\Delta_{ac}^*||_\infty$ in all cases. In addition, $||\Delta_{ac}^*||_\infty$ is relatively large when $\Sigma$ is exchangeable, $s_0 = 15$ and $p = 500$, which explains the lack of performance/undercoverage in this case. We observe that the maximum norms of $\Delta_{ac}^*$ and $\Delta_{in}^*$ generally increase with $s_0$. Overall, the above discussions support our observations in Tables 1-2 of the paper in the sense that the lower the (normalized) remainder term is, the more accurate the coverage is.

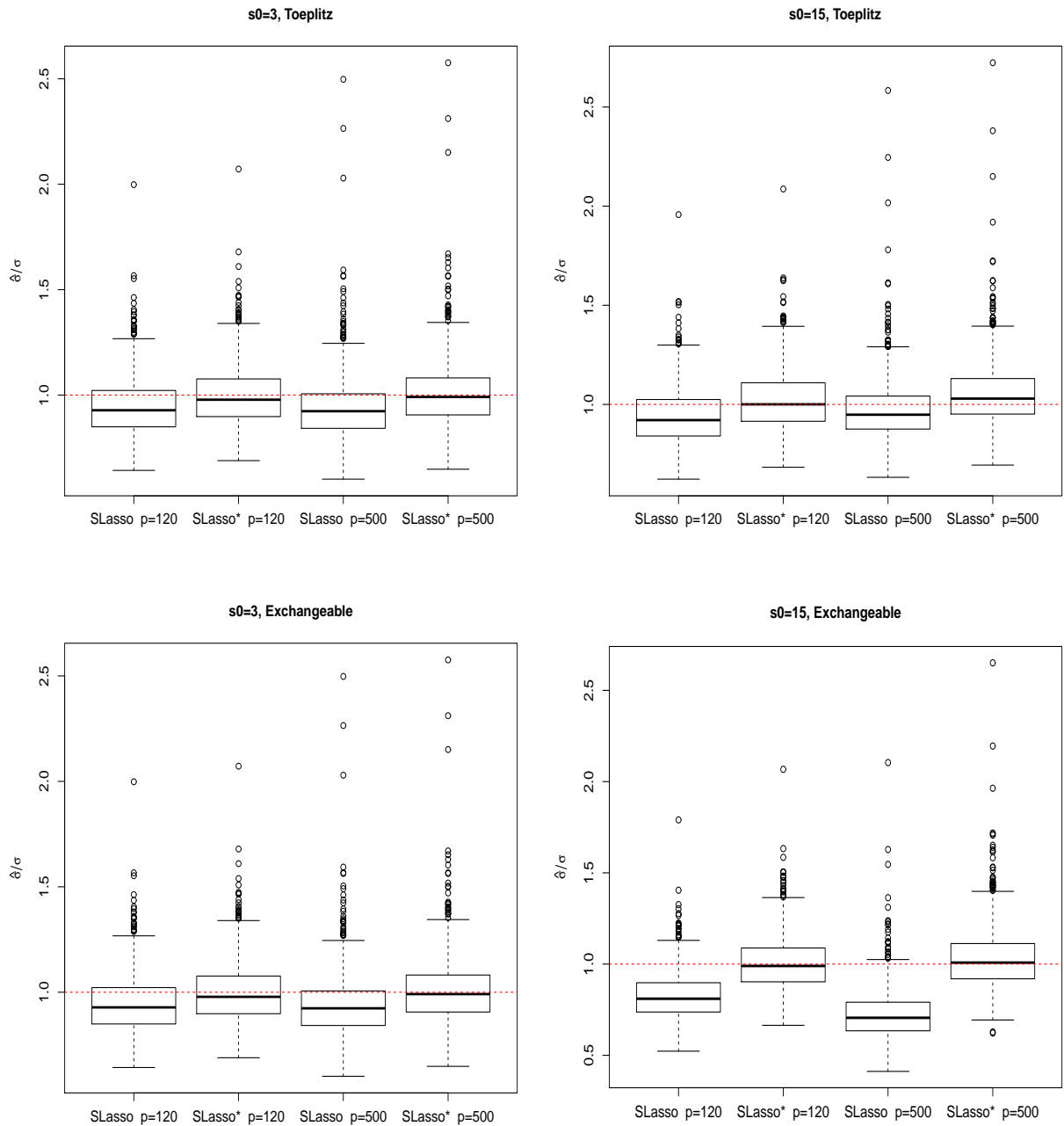

Figure S.1: Boxplots for $\widehat{\sigma}/\sigma$, where $s_0 = 3$ or 15, $\Sigma$ is Toeplitz or exchangeable, and the errors are generated from the studentized $t(4)$ distribution. Here "SLasso" corresponds to the variance estimator delivered by the scaled Lasso and 'SLasso*" corresponds to the modified variance estimator.

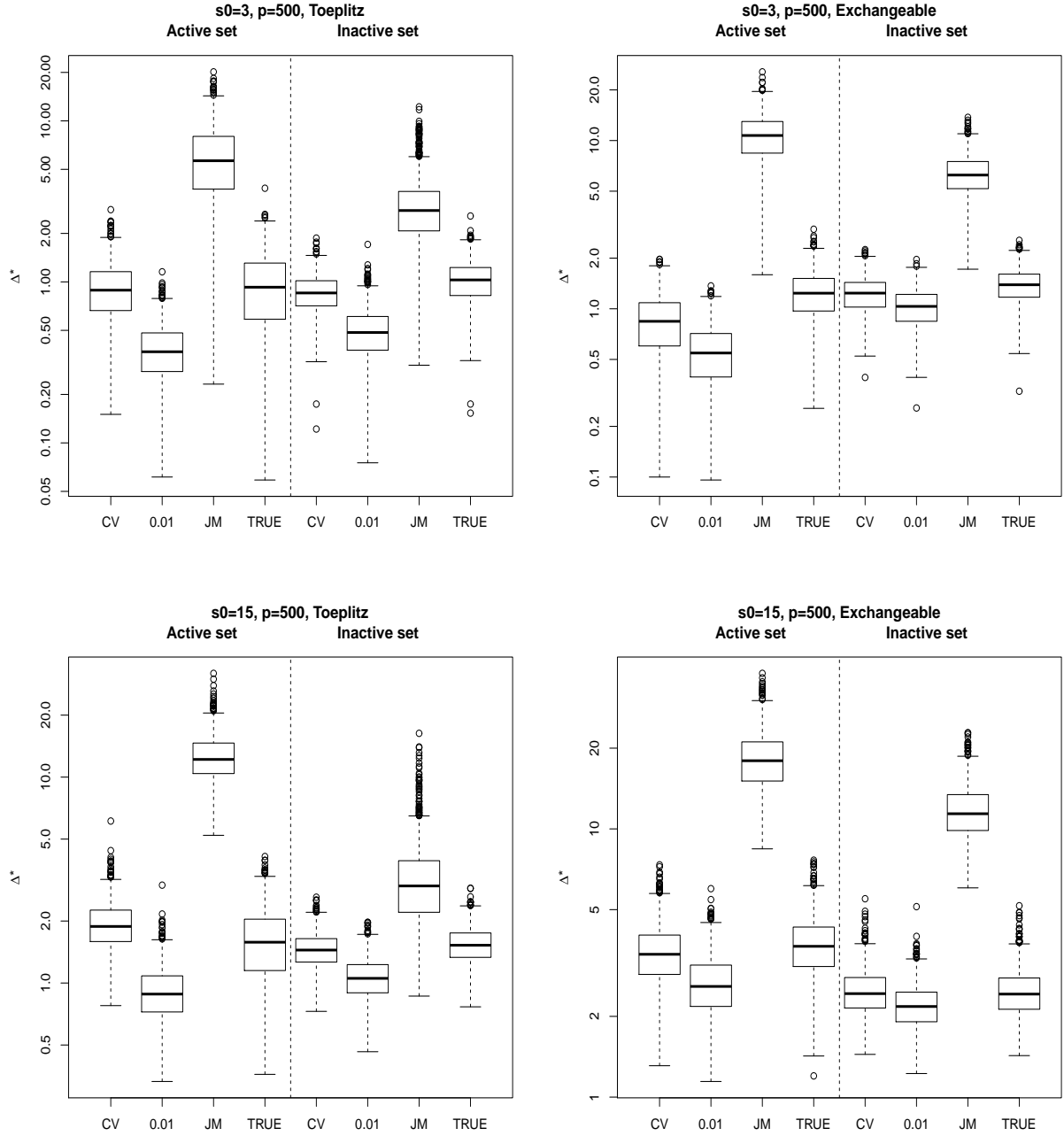

Figure S.2: Boxplots for $||\Delta^*_{\text{ac}}||_\infty$ and $||\Delta^*_{\text{in}}||_\infty$, where $s_0 = 3$ or 15, $p = 500$, $\Sigma$ is Toeplitz or exchangeable, and the errors are $t(4)/\sqrt{2}$. Here "CV", "0.01", "JM" and "TRUE" denote the nodewise Lasso with $\lambda_j$s chosen by 10-fold cross validation and $\lambda_j = 0.01$, the method in Javanmard and Montanari (2014) and the method with the true $\Theta$ respectively. Note that the $y$-axis is plotted on a log scale.